\documentclass[11pt]{article}

\RequirePackage[OT1]{fontenc}

\RequirePackage[aos]{imsart}
\usepackage{amssymb,amsfonts,latexsym,chicago, epsf, here, graphicx,pictexwd}

\startlocaldefs
\newtheorem{theorem}{Theorem}[section]

\newtheorem{lemma}{Lemma}[section]

\newtheorem{prop}{Proposition}[section]

\newtheorem{corollary}{Corollary}[section]
\newtheorem{defn}{Definition}[section]

\newtheorem{remark}{Remark}[section]

\newtheorem{example}{Example}[section]

\newcommand{\GG}{ {\mathbb G}}

\newcommand{\RR}{ {\mathbb R}}

\newcommand{\mycite}[1]{{\small \sc \citeNP{#1}}}
\endlocaldefs

\begin{document}

\begin{frontmatter}

\title{Estimation of a $k-$monotone density, part 1:
characterizations, consistency, 
and minimax lower bounds}
\runtitle{$k-$monotone: characterizations, consistency, lower bounds}


\author{\fnms{Fadoua} \snm{Balabdaoui}\thanksref{t1,t3}\corref{}\ead[label=e1]{fadoua@math.uni-goettingen.de}}
\thankstext{t1}{Research supported in part by National Science Foundation
grant DMS-0203320}
\address{Institute for Mathematical Stochastics\\
Georgia Augusta University Goettingen\\
Maschmuehlenweg 8-10\\
D-37073 Goettingen\\
Germany \\\printead{e1}}
\and
\author{\fnms{Jon A} \snm{Wellner}\thanksref{t2}\ead[label=e2]{jaw@stat.washington.edu}}
\thankstext{t2}{Research supported in part by National Science Foundation grant
DMS-0203320,  NIAID grant 2R01 AI291968-04, and an NWO Grant
to the Vrije Universiteit, Amsterdam}
\thankstext{t3}{Corresponding author}
\address{Department of Statistics\\Box 354322 \\University of Washington\\Seattle, WA 98195-4322 \\
\printead{e2}}

\affiliation{University of G\"ottingen and University of Washington}

\runauthor{Balabdaoui and Wellner}

\begin{abstract}
Shape constrained densities are encountered in many 
nonparametric estimation problems.  The 
classes of monotone or convex (and monotone)
densities can be viewed as special cases of the classes 
of {\sl $k-$monotone} densities.  A density $g$ is said to be 
$k-$monotone if 
$(-1)^l g^{(l)} $ is nonnegative, nonincreasing
and convex for $l=0, \ldots , k-2$ if $k\ge 2$, and $g $ is simply nonincreasing
if $k=1$.  These classes of shaped constrained densities bridge the gap
between the classes of monotone (1-monotone) and 
convex decreasing (2-monotone) densities for which asymptotic 
results are known, and the class of completely monotone ($\infty-$monotone)
densities.   
It is well-known that a density is completely monotone if and only
 if it is a scale mixture
of exponential densities (Bernstein's theorem).  Thus one motivation 
for studying the problem of estimation of a $k-$monotone density
is to try to gain insight into the problem of estimating a completely 
monotone density.

In this series of four papers we consider both (nonparametric) Maximum
Likelihood estimators and Least Squares estimators of a $k-$monotone
estimator.  In this first part (part 1), we prove existence of the estimators
and give characterizations.  We also establish consistency 
properties, and show that the estimators are splines of order $k$ 
(degree $k-1$) with simple knots.  We further provide asymptotic minimax 
risk lower bounds for estimating a $k-$monotone density $g_0 (x_0)$
and its derivatives $g_0^{(j)}(x_0)$, $j=1, \ldots , k-1$,
at a fixed point $x_0$ under the 
assumption that $(-1)^k g_0^{(k)} (x_0) > 0$.

Part 2 of the series gives algorithms for computation of the estimators 
and an application of the methods to earthquake aftershock data.
In part 3 we describe and establish 
existence of the limiting process $H_k$ which governs the asymptotic 
distribution theory modulo a certain conjecture involving a 
Hermite interpolation problem.  In part 4 we give the limiting distribution
theory in terms of $H_k$, again modulo the same Hermite interpolation 
problem.
\end{abstract}

\begin{keyword}[class=AMS]
\kwd[Primary ]{62G05, 60G99}
\kwd[; secondary ]{60G15, 62E20}
\end{keyword}

\begin{keyword}
\kwd{completely monotone}
\kwd{inversion}
\kwd{minimax risk}
\kwd{mixture models}
\kwd{multiply monotone}
\kwd{nonparametric estimation}
\kwd{rates of convergence}
\kwd{shape constraints}
\end{keyword}

\end{frontmatter}

\section{Introduction}
\setcounter{equation}{0}
Shape constrained densities are encountered in many nonparametric
estimation problems.  
Monotone densities arise naturally via connections with renewal 
theory and uniform mixing (see \mycite{vardi:89} and \mycite{wood-sun:93} for 
examples of the former, and \mycite{wood-sun:93} for the latter in an astronomical 
context).   
Convex densities arise in connection with Poisson process models for bird migration
and scale mixtures of triangular densities; see e.g. \mycite{hamp:87}, 
\mycite{anevski:03}, and \mycite{lsm:91}.   

Estimation of monotone densities on the positive half-line 
$\RR^+ = [0,\infty)$ was initiated by \mycite{gren:56} (with 
related work by \mycite{abers:55}, \mycite{brunk:58}, and \mycite{vaneeden:56}, 
\mycite{vaneeden:57}).  Asymptotic theory of the maximum likelihood 
estimators was developed by \mycite{prakra0:69} with later contributions by 
\mycite{gro:85}, \mycite{gro:89}, \mycite{birge:87}, \mycite{birge:89}, and 
\mycite{kimpol:90}.

Estimation of convex densities on $\RR^+$ was apparently initiated by 
\mycite{anevski:94} (see also \mycite{anevski:03}), and was pursued by 
\mycite{wang:94} and \mycite{jongb:95}.  The limit distribution theory for the 
(nonparametric) maximum likelihood estimator and its first derivative at a fixed point 
was obtained by \mycite{gjw:01b}.  

Our goal here (and in the accompanying papers
\mycite{balabwell:04a}, 
\mycite{balabwell:04b}, and \mycite{balabwell:04c})
 is to develop nonparametric estimators and asymptotic theory 
 for the classes of {\sl $k$-monotone densities} on $[0,\infty)$ defined as 
 follows:  $g$ is a $k-$monotone density on $(0,\infty)$ if  
 $g$ is nonnegative and 
$(-1)^l g^{(l)} $ is nonincreasing and convex for $l \in \{ 0, \ldots , k-2\}$ 
for $k\ge 2$, and simply nonnegative and nonincreasing when $k=1$.
As will be shown in section 2, it follows from the results of \mycite{williamson:56}, 
\mycite{levy:62}, and \mycite{gneit:99} that $g$ is a $k-$monotone density if and 
only if it can be represented as a scale mixture of Beta$(1,k)$ densities; i.e. 
with $x_+ \equiv x 1\{ x\ge 0\}$, 
$$
g(x) = \int_0^{\infty} \frac{k}{y^k} (y-x)_+^{k-1} dF(y) 
$$
for some distribution function $F$ on $(0,\infty)$. 
Note that for $k=1$ this recovers the well known fact that monotone densities 
are in a one-to-one correspondence with scale mixtures of uniform densities, 
and, for $k=2$, the corresponding fact frequently used by \mycite{gjw:01b} 
that convex decreasing densities are in a one-to-one correspondence with 
scale mixtures of the triangular, or Beta$(1,2)$, densities.

Our motivation for studying nonparametric estimation in the classes ${\cal D}_k$
has several components:  besides the obvious goal of generalizing the 
existing theory for the $1-$monotone (i.e. monotone) and $2-$monotone (i.e.
convex and decreasing) classes ${\cal D}_1$ and ${\cal D}_2$, these classes play 
an important role in several extensions of  Hampel's bird migration problem 
which are discussed further in \mycite{balabwell:05a}. 
They also provide a potential link to the important limiting case of the 
$k-$monotone classes, namely the class ${\cal D}_{\infty}$ of completely monotone
densities.  Densities $g$ in ${\cal D}_{\infty}$ have the property that 
$(-1)^l g^{(l)} (x) \ge 0$ for all $x \in (0,\infty)$ and $l \in \{ 0, 1, \ldots \}$.
It follows from Bernstein's theorem (see e.g. \mycite{feller:71}, page 439, 
or \mycite{gneit:98}) that $g \in {\cal D}_{\infty}$ if and only if it can be 
represented as a scale mixture of exponential densities; i.e. 
$$
g(x) = \int_0^{\infty} y^{-1} \exp(-x/y) dF(y)
$$
for some distribution function $F$ on $(0,\infty)$.  
Completely monotone densities arise naturally in connection 
with mixtures of Poisson processes and have been 
used in reliability theory (see e.g. 
\mycite{harr-sing:68}, \mycite{dhm:80}, \mycite{hsl:80}), 
and empirical Bayes estimation (see \mycite{robbins:64} and \mycite{robbins:80}).
\mycite{jewell:82} initiated the study of maximum likelihood estimation in the
family ${\cal D}_{\infty}$ and succeeded in showing that the MLE
$\hat{F}_n$ of the mixing distribution function $F$ is unique and 
almost surely weakly consistent.  Although consistency of the MLE follows 
now rather easily from the results of \mycite{pfanz:88} and 
\mycite{vdg:93}, little is known about rates of convergence or asymptotic
distribution theory for either the 
estimator $\hat{g}_n$ of the mixed density $g$ in ${\cal D}_{\infty}$ 
(the ``forward'' or ``direct'' problem) or the estimator $\hat{F}_n$ of 
the mixing distribution function $F$ (the ``inverse'' problem).  
Although our present methods do not yield solutions of these difficult questions, 
the development of methods and theory for general $k-$monotone densites
may throw some light on the issues and problems.  
\medskip

Now we briefly describe the contents of the four related papers 
of which the present manuscript is part 1.

In this paper (part 1), we consider the Maximum Likelihood $\hat{g}_n$ 
and Least Squares $\tilde{g}_n$ estimators of a density $g_0 \in {\cal D}_k$ 
for a fixed integer $k\ge 2$ based on a sample 
$X_1, \ldots , X_n$ i.i.d. with density $g_0$.  
We show that the estimators exist,
provide 
characterizations, and establish consistency of the estimators and their
 derivatives $\hat{g}_n^{(j)}$ and $\tilde{g}_n^{(j)}$ 
for $j \in \{ 1, \ldots , k-1\}$ (uniformly on closed sets bounded away from 
$0$).   
In section 4 we establish asymptotic
minimax lower bounds for estimation of $g_0^{(j)}(x_0)$, $j= 0, \ldots , k-1$ 
under the assumption that $g_0^{(k)}(x_0)$ exists and is non-zero.  
In part 1 we also include statements of known results for estimation of 
a completely monotone density $g_0 \in {\cal D}_{\infty}$ whenever possible. 
One of the remaining open questions concerns existence of the least squares estimator;
see Section 2. 
In section~\ref{NumIllustr} we illustrate both the maximum likelihood and 
least squaqres estimators for $k=3$ and $k=6$ in both the 
direct and inverse problems via artifical data generated from 
a standard exponential distribution.
\medskip

In part 2 (\mycite{balabwell:04a}) we provide algorithms for computation 
of the estimators and for computation of  (approximations to) the limit
process $H_{c,k}$  defined in part 3 (\mycite{balabwell:04b}).  
We call the basic algorithm developed and used 
in part 2 an {\sl iterative $(2k-1)-$spline algorithm} since
it extends the ``cubic spline algorithm'' developed in 
\mycite{gjw:01a} and \mycite{gjw:03}.
Part 3 is devoted to a study of the corresponding canonical Gaussian 
problem and the ``invelope'' ($k$ even) or envelope ($k$ odd) processes
$H_k = \lim_{c\rightarrow \infty} H_{c,k}$ which arise in the solution 
of the Gaussian version of the problem.  
Thus part 3 extends and is analogous to the treatment for 
the case $k=2$ given by \mycite{gjw:01a}.  
Finally, part 4 (\mycite{balabwell:04c}) gives joint asymptotic distribution 
theory at a fixed point $x_0 \in (0,\infty)$ of the vector of centered and 
scaled derivative estimators 
$$
\left ( n^{(k-j)/(2k+1)} ( \overline{g}_n^{(j)} (x_0) - g_0^{(j)} (x_0) ),  \ \ 
j = 0, \ldots , k-1 \right ) \
$$
where $\overline{g}_n $ is either the MLE $\hat{g}_n$ or the LSE $\tilde{g}_n$,
under the assumption that $g_0^{(k)} (x_0)$ exists and is non-zero.   
This yields behavior of the corresponding estimators of the 
mixing distribution $F_0$ at fixed points (the inverse problem) as a corollary.
\medskip

Thus the main outcome of parts 3 and 4 
generalizes the asymptotic distribution theory for estimating a nondecreasing density, and a nondecreasing and convex density 
at a fixed point: 
If $x_0 > 0$ and $g_0$ is a $k$-monotone density defined on $(0,\infty)$
such that $g_0$ is $k$-times differentiable at $x_0$ with $(-1)^k g^{(k)}_0(x_0) > 0 $, and $g^{(k)}_0$ is assumed to be continuous in a neighborhood of $x_0$, then 
our goal in parts 3 and 4 is to show that 
\begin{eqnarray*} \left (
\begin{array}{c}
  n^{\frac{k}{2k+1}}(\bar{g}_{n}(x_{0})-g_0(x_{0}) ) \\
  n^{\frac{k-1}{2k+1}}(\bar{g}^{(1)}_{n}(x_{0})-g^{(1)}_0(x_{0}) ) \\
  \vdots \\
  n^{\frac{1}{2k+1}}(\bar{g}^{(k-1)}_{n}(x_{0})-g^{(k-1)}_0(x_{0}))
\end{array} \right )
\rightarrow_{d}
\left (
\begin{array}{c}
c_{0}(g_0)H^{(k)}_k(0)\\
c_{1}(g_0)H^{(k+1)}_k(0)\\
\vdots \\
c_{k-1}(g_0)H^{(2k-1)}_k(0)\\
\end{array}\right )\\
\end{eqnarray*}
and
\begin{eqnarray*}
  n^{\frac{1}{2k+1}} (\bar{F}_{n}(x_0)- F_0 (x_0))
\rightarrow_{d}\frac{(-1)^k x^{k}_0}{k!}c_{k-1}(g_0)H_k ^{(2k-1)}(0),
\end{eqnarray*}
where $\bar{g}_n$ is either the MLE of LSE, $\bar{F}_n$ is the corresponding 
estimator of the mixing distribution function $F_0$, 
\begin{eqnarray*}
c_{j}(g_0) = \bigg \{\left(g_0(x_0)\right)^{k-j}
    \left(\frac{(-1)^k g^{(k)}_0(x_0)}{k!}\right)^{2j+1} \bigg \}^{\frac{1}{2k+1}},
\end{eqnarray*}
for $j=0, \cdots, k-1$, and $H_k$ is an almost surely uniquely defined stochastic process that is $(2k)$-convex (i.e., $H^{(2k-2)}_k$ exists and convex), 
and stays above (below) the $(k-1)$-fold integral of two-sided Brownian motion plus a polynomial drift of the form $t^{2k}/(2k)!$ if
 $k$ is even (odd). Only a change of scale is 
 necessary to realize that $H_1$ and $H_2$ are very closely related to the greatest convex minorant of $W(t) + t^2, t \in \mathbb{R}$, where $W$ is two-sided Brownian motion, and the ``invelope'', $H$, of
\begin{eqnarray*}
\left \{
\begin{array}{ll}
 \int_0^t W(s) ds  + t^4, \ \textrm{if $t \geq 0$} \\
 \int_t^0 W(s) ds  + t^4, \ \textrm{if $t < 0$}.
\end{array}
\right. 
\end{eqnarray*}
Deriving the rate of convergence of both the estimators $\hat{g}_n$ 
and $\tilde{g}_n$ and their derivatives
$\hat{g}_n^{(j)}$, $\tilde{g}_n^{(j)}$,  $j=1, \cdots, k-1$, 
and proving the existence of the stochastic processes $H_k$ for $k > 2$
 involved in the joint asymptotic distribution still depends on a key conjecture:  
 that the distance between two successive knots of the MLE or LSE that are 
 in the neighborhood of $x_0$ is $O_p (n^{-1/(2k+1)}) $ as the sample size 
 $n \to \infty$, and that distance between two successive points of touch 
 between the $(k-1)$-fold integral of two-sided Brownian motion plus 
 $t^{2k}/(2k)!$ and $H_k$ is $O_p(1)$. Both problems are of the same 
 nature and one can go from the first to the second one via a simple 
 scaling argument.  
 We refer to this common problem as the \textit{gap problem}.
 
 We will show in parts 3 and 4 that 
 the gap problem can be reduced to the solution of a certain 
 problem related to Hermite interpolation.  
 That is, the gap problem has a solution if the following conjecture
 involving Hermite interpolation is true:
 Consider 
 \textrm{Hermite interpolation} (as 
 described for example in 
 \mycite{Nurn:89}, pages 108-109 or 
  \mycite{devore-lorentz:93} pages 161 - 162)
 of some smooth function $f$ via splines of 
 odd-degree. More specifically, if $f$ is some real-valued function in 
 $C^{(j)}[0,1]$ for some $j \geq 2$, $0= y_0 < y_1 < \cdots < y_{2k-4} < y_{2k-3} = 1$ 
 is a given increasing sequence, then the uniquely defined spline $Hf$ of 
 degree $2k-1$ and interior knots $y_1, \ldots, y_{2k-4}$ satisfying the $4k-4$ conditions
\begin{eqnarray*}
(Hf)(y_i)= f(y_i), \ \ \textrm{and} \ \ (Hf)'(y_i) = f'(y_i), \ \  i=0, \ldots, 2k-3,
\end{eqnarray*}
then we conjecture that there exists a constant $c_{k,j}$ depending only on $k$ and $j$ such that, if $j\ge k$,   
\begin{eqnarray*}
\sup_{0 < y_1< \cdots < y_{2k-4} < 1} \Vert f - Hf \Vert_\infty \leq c_{k,j} \ \Vert f^{(j)} \Vert_\infty, 
\end{eqnarray*}
where $\Vert \cdot \Vert_\infty$ is the supremum norm over $[0,1]$. 

This Hermite interpolation problem has apparently not been 
investigated in detail in the 
spline or approximation theory 
literature, and hence an analysis of the corresponding
interpolation error is yet to be developed.   It is, however, precisely the 
interpolation problem involved in understanding our least squares estimators,
both for finite sample sizes and in the limiting Gaussian problem:  as will 
be shown in parts 3 and 4,
the connecting link is the classical theorem of 
\mycite{schoen-whit:53} and its generalization by 
\mycite{karlin-zieg:66}; see \mycite{Nurn:89}, page 109, or 
\mycite{devore-lorentz:93}, page 162.

However, the approximation theory literature has considered
a related conjecture for another Hermite 
problem whose solution is a different odd-degree spline, also 
called a \textit{complete spline}. Given a function $f \in C^{(k-1)}[0,1]$, 
and  an increasing sequence $0= y_0 < y_1 < \cdots < y_{m} < y_{m+1} =1$,  
the complete spline interpolant, \ $Cf$,  of degree $2k-1$ with interior 
knots $y_1, \cdots, y_m$ satisfies the $2k + m$ conditions
\begin{eqnarray*}
\left \{
\begin{array}{ll}
(Cf)(y_i) = f(y_i), \ \ i=1, \cdots, m \\
(Cf)^{(l)}(y_0) = f^{(l)}(y_0), \ (Cf)^{(l)}(y_{m+1}) = f^{(l)}(y_{m+1}), \ \ l=0, \cdots, k-1.
\end{array}
\right.
\end{eqnarray*}
When $f$ is in $C^{(j)}[0,1]$ for $j\ge k$, 
the error in this more usual 
Hermite problem is known to be uniformly bounded independently 
of the location of the knots. 
Proof of this uniform boundedness is due to \mycite{shadrin:92}.
More precisely, the argument follows from his Theorem 6.4, page 94. 
\mycite{deboor:74} had investigated the problem for $j=2k$, and 
conjectured uniformity of the bound for this particular case.  
Furthermore, \mycite{deboor:74}  reduced the problem to a further conjecture:
for any $k>4$, the supremum norm of the $L_2-$spline projector that maps
$C^{(k)}[0,1]$ to the space of splines of degree $k-1$ with knots 
$y_1, \ldots, y_m$ is bounded independently of the location of the knots.
This conjecture remained unsolved for more than 25 years:  \mycite{shadrin:01}
presents a proof thereof.  Thus, there is a closely related interpolation problem
in which the interpolation error does hold uniformly in the knots, and this gives 
some hope that ``uniformity in the knots'' will hold in our problem as well.

 In our Hermite interpolation problem, the spline interpolant matches 
 not only the value of the function at the knots but also the value of its 
 first derivative. So intuitively, one should expect our spline to 
 ``behave better''  than the complete spline, and the 
 interpolation error to be smaller.  On the other hand, our conjecture is 
 supported by numerical evidence for $k=3,4,5,6$.  
 Our computations suggest that for these particular values, 
 $c_{k,j} \le 1/((k-1)! (j-k)!)$.
For further details see \mycite{balabwell:05b}.
\bigskip

\section{The Maximum Likelihood and Least Squares estimators: Existence and characterization}
\setcounter{equation}{0}

\subsection{Mixture representation of a $k$-monotone density}

\mycite{williamson:56} gave the following characterization of a
$k$-monotone function on $(0, \infty)$:
\par\noindent
\begin{theorem} (Williamson, 1956)
A function $g$ is $k$-monotone on $(0,\infty)$ if and only if there exists a nondecreasing function
$\gamma$ bounded at 0 such
that
\begin{eqnarray}\label{KMonochar}
g(x) = \int_0^\infty(1-tx)^{k-1}_+ d\gamma(t), \hspace{0.5cm} x > 0
\end{eqnarray}
where $y_{+} = y 1_{(0,\infty)}(y)$.
\end{theorem}
\medskip
The next theorem gives an inversion formula for the function $\gamma$:

\par\noindent
\begin{theorem} (Williamson, 1956) If $g$ is of the form (\ref{KMonochar}) with
$\gamma(0)= 0$, then at a continuity point $t > 0$, $\gamma$ is given by
\begin{eqnarray*}
\gamma(t)  =   \sum_ {j=0}^{k-1}
\frac{(-1)^{k-l}g^{(j)}(1/u)}{j!}\bigg(\frac{1}{u}\bigg)^{j}.
\end{eqnarray*}

\end{theorem}

\medskip

\par\noindent
For proofs of 
Theorems 2.2.1 and 2.2.2, see \mycite{williamson:56}.  \hfill
$\blacksquare $  \\

 From the characterization given in (\ref{KMonochar}),
we can easily derive another integral representation for $k$-monotone 
functions that are
Lebesgue integrable on $(0,\infty)$; i.e., $\int_{0}^{\infty} g(x)dx < \infty$.
\medskip

\par\noindent 
\begin{lemma}
(Integrable $k-$monotone characterization)
A function
$g$ is an integrable $k$-monotone function if and
only if it is of the form
\begin{eqnarray}
g(x) = \int_{0}^\infty  \frac{k (t-x)_{+}^{k-1}}{t^{k}}dF(t), 
\hspace{0.5cm} x > 0
\label{KMonoIntchar}
\end{eqnarray}
where $F$ is nondecreasing and bounded on $(0,\infty)$.
Thus $g$ is a $k-$monotone density if and only if it is of the 
form (\ref{KMonoIntchar}) for some distribution function $F$ on $(0,\infty)$.
\end{lemma}
\medskip

\par\noindent
\textbf{Proof.} 
This follows from Theorem 5 of \mycite{levy:62} by taking $k = n+1$ and
$f\equiv 0 $ on $(-\infty, 0 \rbrack$.
\hfill $\blacksquare$
\medskip

\par\noindent
\begin{lemma}
\label{Inversion}
(k-monotone inversion formula)
If $F$ in (\ref{KMonoIntchar}) satisfies $\lim_{t\to \infty}
F(t) = \int_{0}^\infty g(x)dx$, then
at a continuity point $t > 0$, $F$
is given by
\begin{eqnarray}
\label{Inverseformula}
\phantom{bla}F(t) = G(t) - tg(t)
   + \cdots + \frac{(-1)^{k-1}}{(k-1)!}t^{k-1}g^{(k-2)}(t) +
      \frac{(-1)^{k}}{k!}t^{k}g^{(k-1)}(t),
\end{eqnarray}
where $G(t) = \int_{0}^t g(x)dx$. \\
\end{lemma}
\medskip

\par\noindent
\textbf{Proof.} 
By the mixture form in (\ref{KMonoIntchar}), we have for all $t > 0$
\begin{eqnarray*}
F(\infty)-F(t)  =  \frac{(-1)^{k}}{k!}\int_{t}^{\infty}x^{k} dg^{(k-1)}(x).
\end{eqnarray*} 
But, for $j=1, \cdots, k$, $t^j G^{(j)}(t) \searrow 0$ as $t \to \infty$. This follows from Lemma 1 in \mycite{williamson:56} applied to the $(k+1)$-monotone function $G(\infty) - G(t)$. Therefore, for $j=1, \cdots, k$, $t^j g^{(j-1)}(t) \searrow 0$ as $t \to \infty$.

Now, using integration by parts, we can write
\begin{eqnarray*}
\lefteqn{F(\infty)-F(t) }\\
& = & \frac{(-1)^{k}}{k!}\bigg \lbrack  x^{k} g^{(k-1)}(x) 
            \bigg \rbrack_{t}^{\infty} + \frac{(-1)^{(k-1)}}{(k-1)!} \int_{t}^{\infty} x^{k-1}g^{(k-1)}(x)dx \\
& = & - \frac{(-1)^{k}}{k!} t^{k}g^{(k-1)}(t) - \frac{(-1)^{k-1}}{(k-1)!} t^{k-1}g^{(k-2)}(t)\\
 \hspace{0.2cm}&& + \hspace{0.08cm} \frac{(-1)^{k-2}}{(k-2)!} \int_{t}^{\infty}x^{k-2}g^{(k-2)}(x)dx \\
& \vdots & \\
& = & - \frac{(-1)^{k}}{k!} t^{k}g^{(k-1)}(t) - \frac{(-1)^{k-1}}{(k-1)!} t^{k-1}g^{(k-2)}(x)+  \cdots - \int_{t} ^{\infty} g(x)dx,
\end{eqnarray*}
Using the fact that $F(\infty) = \int _{0}^{\infty}g(x)dx$, the result follows.
\hfill $\blacksquare$
\medskip

For completeness and for comparison, we also give the corresponding characterization 
and inversion formula in the completely monotone case:
\medskip

\par\noindent 
\begin{lemma}
\label{ComplMonCharacterization}
(Integrable completely monotone characterization)
A function
$g$ is an integrable completely monotone function if and
only if it is of the form
\begin{eqnarray}
g(x) = \int_{0}^\infty  \frac{1}{t} \exp(-x/t) dF(t), 
\hspace{0.5cm} x > 0
\label{ComplMonoIntchar}
\end{eqnarray}
where $F$ is nondecreasing and bounded on $(0,\infty)$.
Thus $g$ is a completely monotone density if and only if it is of the 
form (\ref{ComplMonoIntchar}) for some distribution function $F$ on $(0,\infty)$.
\end{lemma}
\medskip

\par\noindent
\begin{lemma}
\label{Inversion-ComplMon}
(Completely-monotone inversion formula)
If $F$ in (\ref{ComplMonoIntchar}) satisfies $\lim_{t\to \infty}
F(t) = \int_{0}^\infty g(x)dx$, then
at a continuity point $t > 0$, $F$
is given by
\begin{eqnarray}F(t) = \lim_{k\rightarrow \infty}
\sum_{j=0}^{k} \frac{(-1)^j }{j!} (kt)^j G^{(j)} (kt)
\label{ComplMonInverseformula}
\end{eqnarray}
where $G(t) = \int_{0}^t g(x)dx$. \\
\end{lemma}
\medskip

\par\noindent
{\bf Proofs.}  
Lemma~\ref{ComplMonCharacterization} follows from 
 the 
classical result of Bernstein; see 
\mycite{widder:46}, pages 141-163; 
\mycite{feller:71}, page 439; and  
\mycite{gneit:98}.
 Lemma~\ref{Inversion-ComplMon} 
  follows from the development in \mycite{feller:71}, pages 232-233.
  For further details, see \mycite{balabwell:05a}. 
  \hfill $\blacksquare$

\medskip

The characterization in (\ref{KMonoIntchar}) is
more relevant for us since we are dealing with $k$-monotone 
densities. It is easy to see that
if $g$ is a density, and $F$ is chosen to be right-continuous and to satisfy the condition of Lemma \ref{Inversion}, then $F$ is a distribution function.
For $k=1$ ($k=2$), note that the
characterization matches with the well known fact that a density is 
nondecreasing
(nondecreasing and convex) on $(0,\infty)$ if and only if it is a 
mixture of uniform
densities (triangular densities). More generally, the 
characterization establishes a
one-to-one correspondance between the class of $k$-monotone densities 
and the class of scale
mixture of Beta's with parameters $1$ and $k$.  From the inversion
formula in (\ref{Inverseformula}), one can see that a natural estimator for the 
mixing distribution $F$ is
obtained by plugging in an estimator for the density $g$ and it 
becomes clear that the rate of convergence of estimators of
$F$ will be controlled by the corresponding rate of 
convergence for estimators of the highest 
derivative $g^{(k-1)}$ of $g$.
When $ k $ increases
  the densities become smoother, and therefore the
inverse problem of estimating the mixing distribution $F$ becomes harder.

In the next section, we
consider the nonparametric Maximum Likelihood and Least Squares Estimators
of a $k$-monotone density $g_0$.  We show that these estimators
exist and give  characterizations thereof.
In the following, $\mathcal{M}_k$ is the class of all
$k$-monotone functions on $(0,\infty)$,
$\mathcal{D}_k$ is the sub-class of 
$k$-monotone densities on $(0,\infty)$, $X_1, \cdots, X_n$ are
i.i.d. from  $g_0$,  and $\mathbb{G}_n$ is their empirical distribution function, 
$\GG_n (x) = n^{-1} \sum_1^n 1\{ X_i \le x \}$ for $x \ge 0$

\subsection{Maximum likelihood estimation of a $k$-monotone density}

Let
\begin{eqnarray*}
l_{n}(g)= \int_0^\infty \log g(x) \, d\mathbb{G}_n(x) 
\end{eqnarray*}
be the  log-likelihood function (really $n^{-1} $ times the log-likelihood function, but 
we will abuse notation slightly in this same way throughout).   We want to maximize 
$l_n (g)$ over $g \in {\cal D}_k$.  To do this, it is frequently of help to change the optimization 
problem to one over the whole cone ${\cal M}_k \cap L_1 (\lambda)$.  This can be done by 
introducing the ``adjusted likelihood function'' $\psi_n (g)$ defined as follows:
\begin{eqnarray*}
\psi_{n}(g)= \int_0^\infty \log g(x) \, d\mathbb{G}_n(x)  - 
\int_0^\infty g(x) dx,
\end{eqnarray*}
for $g \in \mathcal{M}_k \cap L_1(\lambda)$.  Then, as in GJW (2001a), Lemma 2.3,
page 1661, the maximum likelihood estimator $\hat{g}_n$ also maximizes $\psi_n (g)$ 
over  ${\cal M}_k \cap L_1 (\lambda) $ 
 
 Using the integral representations established in 
the previous subsection,
$\psi_n$ can also be rewritten as
\begin{eqnarray*}
\psi_{n}(F)= 
\left \{ 
\begin{array}{l }
\int_0^\infty \log
\left( \int_0^\infty \frac{k(t-x)^{k-1}_+}{t^k} dF(t)\right) 
d\mathbb{G}_n(x)  -
\int_0^\infty \int_0^\infty \frac{k (t-x)^{k-1}_+}{t^k} dF(t) dx,\\
\int_0^{\infty} \log \left ( \int_0^{\infty} \frac{1}{t} \exp(-x/t) dF(t) \right )
d \mathbb{G}_n (x) \\
\qquad  - \int_0^\infty \int_0^\infty \frac{1}{t} \exp(-x/t) dF(t) dx,
\end{array} \right .
\end{eqnarray*}
where $F$ is bounded and nondecreasing.
\medskip

\par\noindent
\begin{lemma}
The maximum likelihood estimator $\hat{g}_{n,k}$ in the classes 
${\cal D}_k $, $k \in \{1, 2, \ldots , \infty \}$ exists.
Furthermore, $\hat{g}_{n,k}$ is the maximizer of $\psi_n$ over ${\cal M}_k \cap L_1 (\lambda)$.
Moreover, for $k \in \{1, 2, \ldots \}$ 
the density $\hat{g}_{n,k}$ is of the form
\begin{eqnarray*}
\hat{g}_{n,k} (x) =\hat{ w} _{1} \frac{k(\hat{a}_{1}-x)^{k-1}_{+}}{\hat{a}_{1}^{k}}
     + \cdots + \hat{w}_{m}\frac{k(\hat{a}_{m}-x)^{k-1}_{+}}{\hat{a}_{m}^{k}},
\end{eqnarray*}
for some $m= \hat{m}_k$, 
while for $k=\infty$, $\hat{g}_{n,\infty} $ is of the form
\begin{eqnarray*}
\hat{g}_{n,\infty} (x) =\frac{\hat{ w} _{1}}{ \hat{a}_1} \exp(-x/\hat{a}_1 ) 
     + \cdots + \frac{\hat{w}_{m}}{ \hat{a}_m }\exp(-x/\hat{a}_m ) 
\end{eqnarray*}
for some $m= \hat{m}_{\infty}$ 
where $\hat{w}_{1},\cdots, \hat{w}_{m}$ and $\hat{a} _{1},\cdots,\hat{a}_{m}$ are
respectively the weights and the support points of the
maximizing mixing distribution
$\hat{F}_{n,k} $.
\end{lemma}

\medskip

\par\noindent
\textbf{Proof.} 
First, we prove that there
exists a density $\hat{g}_{n}$ that maximizes the \lq\lq usual\rq\rq \hspace{0.01cm} log-likelihood 
$l_{n}= \int_0^\infty \log g(x) d\mathbb{G}_n(x)$ over the class $\mathcal{D}_k$ with 
$k$ finite.   For $g$ in
$\mathcal{D}_{k}$, let $F$ be the distribution function such that
\begin{eqnarray*}
g(x) = \int _{0} ^{\infty} \frac{k(y-x)^{k-1}_{+}}{y^{k}} dF(y).
\end{eqnarray*}
The unicomponent likelihood curve $\Gamma$ as defined by
\mycite{lindsay:83a} (see also \mycite{lindsay:95})
 is then
\begin{eqnarray*}
\Gamma = \bigg \{ \bigg (\frac{k(y-X_{1})^{k-1}_{+}}{y^{k}},
\frac{k(y-X_{2})^{k-1}_{+}}{y^{k}},\cdots,\frac{k(y-X_{n})^{k-1}_{+}}{y^{k}} 
\bigg) : \
\ y \in \lbrack 0,\infty ) \bigg \}.
\end{eqnarray*}
It is easy to see that $\Gamma $ is bounded (notice that the
$i$-th component is equal to 0 whenever $y < X_{i}$).
Also, $\Gamma$ is closed. By Theorems 18
and 22 of \mycite{lindsay:95}, there exists a unique maximizer of
$l_{n}$ and the maximum is
achieved by a discrete distribution function that has at most $n$ 
support points.

Now, let $g$ be a $k$-monotone function in $\mathcal{M}_k \cap L_1(\lambda)$
and let $\int _{0}^{\infty}g(x)dx = c$ so that
$g/c \in \mathcal{D}_{k}$. We have

\begin{eqnarray*}
\psi _{n}(g) - \psi _{n}(\hat{g}_{n})
& = & \int _{0}^{\infty} \log \bigg(\frac{g(x)}{c} \bigg)d\mathbb{G}_{n}(x)
      + \log(c) - c +1 \\
&& \qquad       - \ \int _{0}^{\infty} \log \big(\hat{g}_{n}(x) \big)d\mathbb{G}_{n}(x) \\
& \leq & \int_{0}^{\infty} \log \bigg(\frac{g(x)}{c}
          \bigg)d\mathbb{G}_{n}(x)  - \int_{0}^{\infty} \log
           \big(\hat{g}_{n}(x) \big)d\mathbb{G}_{n}(x) \\
& \leq & 0
\end{eqnarray*}
since $\log(c) \leq c-1$. Thus $\psi_n$ is maximized over $\mathcal{M}_k \cap L_1(\lambda)$ by $\hat{g}_n \in \mathcal{D}_k$.   

In the case $k=\infty$, the assertions of the lemma are proved by \mycite{jewell:82}.
\hfill $\blacksquare$

\bigskip

The following lemma gives a necessary and sufficient condition for a point $t$ 
to be in the support of the maximizing distribution function $\hat{F}_{n,k}$.  
For $k \in \{ 3, \ldots \}$ it generalizes lemma 2.4, 
page 1662, \mycite{gjw:01b}.
\medskip

\par\noindent 
\begin{lemma}
\label{CharaSupp}
Let $X_{1},\cdots,X_{n}$ be i.i.d. random variables from the
true density $g_0$, and let $\hat{F}_{n,k}$ and
$\hat{g}_{n,k}$ be the MLE of the mixing and mixed
distribution respectively. Then, for $k \in \{1, 2, \ldots \}$, 
\begin{eqnarray}\label{CharacMLE}
\hat{H}_{n,k} (t) \equiv \GG_n \left ( 
\frac{k(t-X)^{k-1}_{+}/t^{k}}{\hat{g}_{n,k}(X)} \right ) \leq 1,
\end{eqnarray}
with equality if and only if
$t \in \textrm{supp}(\hat{F}_{n,k})= \{\hat{a}_{1},\cdots,\hat{a}_{m} \}$.
In the case $k=\infty$ 
\begin{eqnarray}\label{CharacMLECM}
\hat{H}_{n,\infty} (t) \equiv \GG_n \left ( 
\frac{\exp(-X/t)}{t \hat{g}_{n,\infty} (X)} \right ) \leq 1, \qquad \mbox{for all} \ \ t> 0 
\end{eqnarray}
with equality if and only if
$t \in \textrm{supp}(\hat{F}_{n,\infty})
= \{\hat{a}_{1},\cdots,\hat{a}_{m} \}$.
\end{lemma}
\medskip

\begin{remark}
By factoring out $t^{k-1}$ and replacing $t$ by $kv$ (say), 
it becomes clear that the function $\hat{H}_{n,\infty}$ on the 
right side of  (\ref{CharacMLECM}) is a natural limiting version 
as $k\rightarrow \infty$ of the functions 
$\hat{H}_{n,k}$ on the right side of (\ref{CharacMLE}).
\end{remark}

\par\noindent
\textbf{Proof.}
Since $\hat{F}_{n}$ maximizes the log-likelihood
\begin{eqnarray*}
l_{n}(F) = \frac{1}{n} \sum _{j=1}^{n}
\log \bigg( \int _{0}^{\infty} \frac{k(y-X_{j})^{k-1}_{+}}{y^{k}}dF(y)\bigg),
\end{eqnarray*}
it follows that for all $t > 0 $
\begin{eqnarray*}
\lim_{\epsilon \searrow 0} \frac{l_{n}((1-\epsilon)\hat{F}_{n}
     + \epsilon \delta _{t}) - l_{n}(\hat{F}_{n})}{\epsilon} \leq 0.
\end{eqnarray*}
This yields
\begin{eqnarray*}
\frac{1}{n} \sum _{j=1} ^{n}
\frac{k(t-X_{j})^{k-1}_{+}/t^{k}- 
\hat{g}_{n}(X_{j})}{\hat{g}_{n}(X_{j})} \leq 0
\end{eqnarray*}
or
\begin{eqnarray}\label{MLEineq}
\frac{1}{n}\sum _{j=1}^{n}
   \frac{k(t-X_{j})^{k-1}_{+}/t^{k}}{\hat{g}_{n}(X_{j})} \leq 1.
\end{eqnarray}
 Now, let $ M_{n} $ be the set defined by
\begin{eqnarray*}
M_{n} = \bigg \{t > 0 \hspace{0.2cm} : \hspace{0.2cm}
\frac{1}{n}\sum _{j=1}^{n} 
\frac{k(t-X_{j})^{k-1}_{+}/t^{k}}{\hat{g}_{n}(X_{j})}
= 1 \bigg \}.
\end{eqnarray*}
We will prove now that $M_{n} = \textrm{supp}(\hat{F}_{n})$.
We write $P_{\hat{F}_{n}}$ for the probability measure associated with 
$\hat{F}_{n}$.
Integrating the left hand side of (\ref{MLEineq}) with respect to $\hat{F}_{n}$, we have
\begin{eqnarray*}
\frac{1}{n}\sum _{j=1}^{n} \frac{\int _{0}^{\infty}
     \bigg( k (t-X_{j})^{k-1}_{+} / t^{k} \bigg) 
d\hat{F}_{n}(t)}{\hat{g}_{n}(X_{j})}
& =  &
\frac{1}{n}\sum _{j=1}^{n} \frac{\hat{g}_{n}(X_{j})}{\hat{g}_{n}(X_{j})} =  1.
\end{eqnarray*}
 But, using the definition of $M_n$, we can write,
\begin{eqnarray*}
1 & = & \frac{1}{n}\sum _{j=1}^{n}
          \frac{\int _{0}^{\infty}\bigg( k (t-X_{j})^{k-1}_{+} / t^{k} \bigg)
           d\hat{F}_{n}(t)}{\hat{g}_{n}(X_{j})} \label{MLEeq} \\
& = &  P_{\hat{F}_{n}}(M_{n}) + \frac{1}{n}\sum _{j=1}^{n}
         \int _{\mathbb{R}^{+} \setminus M_{n}}
         \frac{\bigg( k(t-X_{j})^{k-1}_{+}/t^{k} \bigg)}
     {\hat{g}_{n}(X_{j})} d\hat{F}_{n}(t),
\end{eqnarray*}
and so
\begin{eqnarray*}
P_{\hat{F}_n}(\mathbb{R}^{+} \setminus M_n)
& = & \int_{\mathbb{R}^{+}\setminus M_n}  \frac{1}{n}\sum _{j=1}^{n} \frac{\bigg(
         k(t-X_{j})^{k-1}_{+}/t^{k} \bigg)}{\hat{g}_{n}(X_{j})} 
d\hat{F}_{n}(t) \\
& < & P_{\hat{F}_n}(\mathbb{R}^{+} \setminus M_n), \hspace{0.2cm} 
\textrm{if} \hspace{0.1cm}
         P_{\hat{F}_n}(\mathbb{R}^{+} \setminus M_n) > 0.
\end{eqnarray*}
This is a contradiction and we conclude that
$P_{\hat{F}_n}(\mathbb{R}^{+} \setminus M_n) =0$.

The proof of the result for $k=\infty$ is given by \mycite{jewell:82}, 
page 481.
\hfill $\blacksquare$
\vspace{0.3cm}

\subsection{The Least Squares estimator of a $k$-monotone density}

The least squares criterion is
\begin{equation}
Q_n (g) =  \frac{1}{2} \int_0^{\infty} g^2 (x) dx
- \int_0^{\infty} g(x) d\GG_n (x) \, .
\end{equation}
We want to minimize this over
$g \in {\cal D}_k \cap L_2(\lambda) $, the subset of square
integrable $k-$monotone functions.   
Although existence of a minimizer of $Q_n$ over 
${\cal D}_k \cap L_2 (\lambda)$ is quite easily established, 
the minimizer has a somewhat complicated characterization 
due to the density constraint $\int_0^{\infty} g(x) dx =1$.  
Therefore we 
will actually consider the alternative optimization 
problem of minimizing $Q_n (g)$ over
${\cal M}_k \cap L_2(\lambda) $.  
In this optimization problem existence requires more work, 
but the resulting characterization of the estimator is 
considerably simpler.   Further we will show that even though the
resulting estimator does not necessarily have total mass one, 
it does have total mass converging almost surely to one and it 
consistently estimates $g_0 \in {\cal D}_k$.

Using arguments similar to those in the proof of Theorem 1 in \mycite{williamson:56}, one can show that $g \in {\cal M}_k$ if and only if
$$
g (x) = \int_0^{\infty} (t-x)_+^{k-1} d\mu (t)
$$
for a positive measure $\mu$ on $(0,\infty)$.
Thus we can rewrite the criterion in terms of the
corresponding measures $\mu$:
by Fubini's theorem
\begin{eqnarray*}
\int_0^{\infty} g^2 (x) dx
& = & \int_0^{\infty} \int_0^{\infty} r_k (t,t') d\mu (t) d\mu (t')
\end{eqnarray*}
where
$$
r_k (t,t')
   \equiv  \int_0^{\infty}  (t-x)_+^{k-1} (t'-x)_+^{k-1} dx
= \int_0^{t\wedge t'} (t-x)^{k-1} (t'-x)^{k-1} dx \, ,
$$
and
\begin{eqnarray*}
\int_0^{\infty} g(x) d\GG_n (x)
& = & \int_0^{\infty} \int_0^{\infty} (t-x)_+^{k-1} d\mu (t)  d\GG_n (x) 
      = \int_0^{\infty} s_{n,k} (t) d \mu (t) 
\end{eqnarray*}
where 
$$
s_{n,k} (t) \equiv \GG_n ( (t-X)_+^{k-1}) \, .
$$
Hence it follows that, with $g = g_{\mu}$
\begin{eqnarray*}
Q_n (g)
  =  \frac{1}{2}
          \int_0^{\infty} \int_0^{\infty} r_k (t,t') d\mu (t) d\mu (t')
           - \int_0^{\infty} s_{n,k} (t)
           d \mu (t)
  \equiv  \Phi_n (\mu)
\end{eqnarray*}
Now we want to minimize $\Phi_n$ over the set ${\cal X} $
of all non-negative
measures $\mu$ on $R^+$.
Since $\Phi_n$ is convex and
can be restricted to a subset ${\cal C}$
of ${\cal X}$ on which it is lower semicontinuous,
a solution exists and is unique.
\medskip

\par\noindent 
\begin{prop}
\label{ExistenceLSE}  
The problem of minimizing
$\Phi_n(\mu)$ over all non-negative measures $\mu$
has a unique solution $\tilde{\mu}$.
\end{prop}
\medskip

\par\noindent
{\bf Proof.}
Existence follows from \mycite{zeid:85}, Theorem 38.B, page 152.
Here we verify the hypotheses of that theorem.

We identity
$X$ of Zeidler's theorem with the space ${\cal X}$ of
nonnegative measures on $[0,\infty)$,
and we show that we can take $M$ of Zeidler's theorem
to be
$$
{\cal C} \equiv
\{ \mu \in {\cal X} : \ \mu(t,\infty) \le Dt^{-(k-1/2)} \}
$$
for some constant $D< \infty$.

First, we can, without loss, restrict
the minimization to the space of non-negative measures
on $[X_{(1)},\infty)$ where $X_{(1)} >0$ is the
first order statistic of the data.  To see this,
note that we can decompose
any measure $\mu$ as $\mu= \mu_1 + \mu_2$ where $\mu_1 $
is concentrated on $[0,X_{(1)})$ and $\mu_2$ is concentrated
on $[X_{(1)}, \infty)$.  Since the second term of $\Phi_n$
is zero for $\mu_1$, the contribution of the $\mu_1$ component
to $\Phi_n (\mu)$ is always non-negative, so we make $\inf \Phi_n (\mu)$
no larger by restricting to measures on $[X_{(1)}, \infty)$.

We can restrict further to measures $\mu$
with $\int_0^{\infty} t^{k-1} d \mu (t) \le D$
for some finite $D = D_{\omega}$.
To show this, we first give a lower bound for
$r_k (s,t)$.

For $s,t \ge t_0>0$ we have
\begin{eqnarray}
r_k (s,t) \ge \frac{(1-e^{-v_0})t_0}{2k} s^{k-1} t^{k-1}
\label{LowerBoundForRk}
\end{eqnarray}
where $v_0 \approx 1.59$.
To prove (\ref{LowerBoundForRk}) we will use the inequality
\begin{equation}
(1 - v/k)^{k-1} \ge e^{-v} , \qquad 0 \le v \le v_0, \ \ k \ge 2\, .
\label{LowerBoundVersionExponentialApprox}
\end{equation}
(This inequality holds by straightforward computation; see \mycite{hallwell:79}, especially
their Proposition 2.)
Thus we compute
\begin{eqnarray*}
r_k (s,t)
& = & \int_0^{\infty} (s-x)_+^{k-1} (t-x)_+^{k-1} \, dx \\
& = & s^{k-1} t^{k-1} \int_0^{\infty} (1-x/s)_+^{k-1} (1-x/t)_+^{k-1} dx \\
& = & \frac{1}{k} s^{k-1} t^{k-1} \int_0^{\infty}
        \left (1- \frac{y}{sk} \right )_+^{k-1}
        \left (1- \frac{y}{t k} \right )_+^{k-1} \, dy\\
& \ge & \frac{1}{k} s^{k-1} t^{k-1} \int_0^{v_0 (t\wedge s)}
            e^{-y/s} e^{-y/t} dy \\
& = & \frac{1}{k} s^{k-1} t^{k-1} \int_0^{v_0 (t\wedge s)}
            e^{-c y}  dy,  \qquad c \equiv 1/s + 1/t  \\
& = & \frac{1}{k} s^{k-1} t^{k-1} \frac{1}{c} \int_0^{v_0 (t\wedge s)}
           c e^{-c y}  dy,   \\
& = & \frac{1}{k} s^{k-1} t^{k-1} \frac{1}{c}
            \left ( 1 - \exp(-c (t\wedge s) v_0 ) \right ) \\
& \ge & \frac{1}{k} s^{k-1} t^{k-1} \frac{1}{c}
            \left ( 1 - \exp(- v_0 ) \right )
\end{eqnarray*}
since
$$
c(s\wedge t ) = \frac{s+t}{st} ( s \wedge t)
= \left \{
\begin{array}{l l}
       (t+s) /t, & s \le t \\
       (t+s) /s, & s \ge t
\end{array} \right \} \ge 1 \, .
$$
But we also have
$$
\frac{1}{c} = \frac{1}{(1/s)+(1/t)} = \frac{st}{s+t} \ge \frac{1}{2}
s \wedge t \ge \frac{1}{2} t_0
$$
for $s,t \ge t_0$,
so we conclude that  (\ref{LowerBoundForRk}) holds.

 From the inequality (\ref{LowerBoundForRk})
  we conclude that for measures $\mu$ concentrated
on $[X_{(1)}, \infty)$ we have
$$
\int \!\!\! \int r_k (s,t) d\mu (s) d \mu (t)
\ge
\frac{(1-e^{-v_0})X_{(1)}}{2k} \left (
     \int_0^{\infty} t^{k-1} d \mu (t)  \right )^2 \, .
$$
On the other hand,
\begin{eqnarray*}
\int_0^{\infty} s_{n,k} (t) d\mu (t)
\le \int_0^{\infty} t^{k-1} d \mu (t) \, .
\end{eqnarray*}
Combining these two inequalities
it follows that for any measure $\mu$
concentrated on $[X_{(1)}, \infty)$
we have
\begin{eqnarray*}
\Phi_n(\mu)
& = & \frac{1}{2} \int \!\!\!\int r_k (t,s) d\mu (t ) d \mu (s)  -
         \int_0^{\infty} s_{n,k} (t) d \mu (t) \\
& \ge &  \frac{(1-e^{-v_0})X_{(1)}}{4k} \left (
     \int_0^{\infty} t^{k-1} d \mu (t)  \right )^2
          - \int_0^{\infty} t^{k-1} d \mu (t)  \\
& \equiv  & A m_{k-1}^2 - m_{k-1} \, .
\end{eqnarray*}
This lower bound is strictly positive
if
$$
m_{k-1} > 1/A = \frac{4k}{(1-e^{-v_0})X_{(1)} }\, .
$$
But for such measures $\mu$ we can make  $\Phi$ smaller by
taking the zero measure.  Thus we may restrict the minimization
problem to the collection of measures $\mu$ satisfying
\begin{equation}
m_{k-1} \le 1/A \, .
\label{BoundForkMinusFirstMoment}
\end{equation}
Now we decompose any measure $\mu$ on
$[X_{(1)}, \infty )$ as $\mu = \mu_1 + \mu_2$
where $\mu_1$ is concentrated on $[X_{(1)}, MX_{(n)}]$
and $\mu_2$ is concentrated on $(MX_{(n)}, \infty)$ for some
(large) $M>0$.
Then it follows that
\begin{eqnarray*}
\Phi_n(\mu)
& \ge & \frac{1}{2} \int \!\!\!\int r_k (t,s) d\mu_2 (t ) d \mu_2 (s)  -
         \int_0^{\infty} t^{k-1} d \mu (t) \\
& \ge & \frac{(1-e^{v_0}) M X_{(n)}}{4k} (M X_{(n)})^{2k-2}
         \mu(MX_{(n)} , \infty)^2 - 1/A \\
& \equiv & B \mu (MX_{(n)} , \infty)^2 - 1/A  >  0
\end{eqnarray*}
if
$$
\mu (MX_{(n)}, \infty)^2 > \frac{1}{AB}
= \frac{4k}{(1-e^{-v_0}) X_{(1)}} \frac{4k}{(1-e^{-v_0})
(MX_{(n)})^{2k-1}} \, ,
$$
and hence we can restrict to measures $\mu$ with
$$
\mu (MX_{(n)}, \infty) \le
\frac{4k}{(1-e^{-v_0}) X_{(1)}^{1/2} X_{(n)}^{k-1/2} }
\frac{1}{M^{k-1/2}} \,
$$
for every $M\ge 1$.
But this implies that $\mu$ satisfies
$$
\int_0^{\infty} t^{k-3/4} d \mu (t) \le D
$$
for some $0< D = D_{\omega} < \infty$, and this implies that
$t^{k-1}$ is uniformly integrable over $\mu \in {\cal C}$.
Alternatively,
for $\lambda \ge 1$ we have
\begin{eqnarray*}
\int_{t > \lambda} t^{k-1} d \mu (t)
& = & \lambda^{k-1} \mu(\lambda , \infty) + (k-1) \int_\lambda ^{\infty}
         s^{k-2} \mu (s,\infty) ds \\
& \le &  \lambda^{k-1} \frac{K}{\lambda^{k-1/2}} +
          (k-1) \int_{\lambda}^{\infty} s^{k-2} K s^{-(k-1/2)} ds \\
& = &  K \lambda^{-1/2} + (k-1)  K \int_{\lambda}^{\infty} s^{-3/2} ds \\
& \le & K \lambda^{-1/2} + (k-1) 2K \lambda ^{-1/2} \\
& \rightarrow & 0 \qquad \mbox{as} \ \ \lambda \rightarrow \infty
\end{eqnarray*}
uniformly in $\mu \in {\cal C}$.

This implies that
   for $\{ \mu_m\} \subset
{\cal C}$ satisfying $\mu_m \Rightarrow \mu_0$ we have
$$
\limsup \int_0^{\infty} s_{n,k} (t) d \mu_m (t)
\le \int_0^{\infty} s_{n,k} (t) d\mu_0 (t) \, ,
$$
and hence $\Phi$ is lower-semicontinuous on ${\cal C}$:
$$
\liminf_{m\rightarrow \infty}
\Phi_n(\mu_m) \ge \Phi (\mu_0) \, .
$$
Since $\Phi_n$ is lower semi-compact (i.e.
the sets ${\cal C}_r \equiv \{ \mu \in {\cal C} : \ \Phi_n(\mu) \le r \}$
are compact for $r \in \RR$), the existence of
a minimum follows
from \mycite{zeid:85}, Theorem 38.B, page 152.
Uniqueness follows from the strict convexity of $\Phi_n$.
   \hfill $\blacksquare$

\medskip

The following proposition characterizes the least squares estimators.
\medskip

\par\noindent 
\begin{prop}
\label{characterizationLSE}
For $k \in \{ 1, 2, \ldots  \}$ 
define $\mathbb{Y}_{n,k}$ and $\tilde{H}_{n,k}$ respectively by
\begin{eqnarray*}
\mathbb{Y}_{n,k}(t) = 
\int _{0} ^{t} \int _{0} ^{t_{k-1}}\cdots
   \int _{0} ^{t_{2}} \mathbb{G}_{n}(t_{1}) dt_{1}dt_{2}\cdots 
    dt_{k-1}, & \qquad   x \geq 0, 
\end{eqnarray*}
and
\begin{eqnarray*}
\tilde{H}_{n,k}(t)
= 
\int _{0} ^{t} \int _{0} ^{t_{k}}\cdots \int _{0} ^{t_{2}} 
\tilde{g}_{n}(t_{1})
      dt_{1}dt_{2}\cdots dt_{k}, & \qquad  x \geq 0\, .
\end{eqnarray*}
Then $\tilde{g}_{n,k}$ is the LS estimator over
$\mathcal{M}_{k} \cap L_2(\lambda)$ if and only if the following conditions
are satisfied for $\tilde{g}_{n,k}$
and $\tilde{H}_{n,k}$:
\begin{eqnarray}\label{Fenchel}
\left \{ 
\begin{array}{lll}
\tilde{H}_{n,k}(t)  \geq \mathbb{Y}_{n,k}(t), \qquad  \textrm{for} \ \  t \geq 0, \ \ \textrm{and} \\
\tilde{H}_{n,k}(t)  = \mathbb{Y}_{n,k}(t),  \qquad  \mbox{for} \ \ t \in \mbox{supp} \{ 
\tilde{F}_{n,k} \} \, .
\end{array}
\right.
\end{eqnarray}
\end{prop}
\medskip

\par \noindent 
\begin{remark} 
Note that 
for $k \in \{1, 2, \ldots \}$ the processes 
$\mathbb{Y}_{n,k}$ and $\tilde{H}_{n,k}$ can be written in the more compact forms
\begin{eqnarray*}
\mathbb{Y}_{n,k}(t) = \int_0^t \frac{(t-x)^{k-1}}{(k-1)!}d\mathbb{G}_n(x)
\end{eqnarray*}
and
\begin{eqnarray*}
\tilde{H}_{n,k}(t) = \int_0^t \frac{(t-x)^{k-1}}{(k-1)!}\tilde{g}_n(x)dx.
\end{eqnarray*}
\end{remark}
\medskip

\par\noindent
\textbf{Proof.}
Let $\tilde{g}_{n} \in \mathcal{M}_{k} \cap L_2(\lambda)$ satisfy
(\ref{Fenchel}), and let $g$ be an arbitrary function in 
$\mathcal{M}_{k} \cap L_2(\lambda)$. Then
\begin{eqnarray*}
Q_{n}(g)-Q_{n}(\tilde{g}_{n})
& = & \frac{1}{2}\int g^{2}(x)dx-\frac{1}{2}\int \tilde{g}_{n}^{2}(x)dx\\
&& \qquad         - \ \int g(x)
      d\mathbb{G}_{n}(x) + \int \tilde{g}_{n}(x) d\mathbb{G}_{n}(x).
\end{eqnarray*}
Now, using integration by parts
\begin{eqnarray*}
\lefteqn{ \int _{0} ^{\infty} (g(x) - \tilde{g}_{n}(x)) d\mathbb{G}_{n}(x)}\\
& = & - \int _{0} ^{\infty} \mathbb{G}_{n}(x) (g'(x) - 
\tilde{g}'_{n}(x))  dx \\
& = & \int _{0} ^{\infty}  \bigg(\int _{0}^{x} \mathbb{G}_{n}(y)dy \bigg)
       (g''(x) - \tilde{g}''_{n}(x)) dx \\
&\vdots & \\
& = & (-1)^{k} \int _{0} ^{\infty} \mathbb{Y}_{n}(x) (dg^{(k-1)}(x) - 
d\tilde{g}^{(k-1)}_{n}(x))  ,
\end{eqnarray*}
and
\begin{eqnarray*}
\lefteqn{\int _{0}^{\infty} (g^{2}(x)-\tilde{g}_{n}^{2}(x)) dx }\\
&= & \int _{0}^{\infty} (g(x)+\tilde{g}_{n}(x))(g(x)- \tilde{g}_{n}(x))dx \\
&=& - \int _{0} ^{\infty} \bigg(\int _{0}^{x} g(y)dy
       +  \int _{0}^{x} \tilde{g}_{n}(y)dy \bigg) (g'(x)- 
\tilde{g}'_{n}(x))dx \\
& \vdots & \\
&=& (-1)^{k} \int _{0} ^{\infty} (G_{k}(x) + \tilde{H}_{n}(x))
       (dg^{(k-1)}(x)-d\tilde{g}^{(k-1)}_{n}(x)),
\end{eqnarray*}
where $G_{k}$ is the $k$-th order integral of $g$. Hence,
\medskip
\begin{eqnarray*}
Q_{n}(g)-Q_{n}(\tilde{g}_{n})
&=& \frac{1}{2} (-1)^{k} \int _{0}^{\infty}
      (G_{k}(x)+\tilde{H}_{n}(x))(dg^{(k-1)}(x)-d\tilde{g}^{(k-1)}_{n}(x)) \\
&&- \ (-1)^{k} \int_{0}^{\infty} \mathbb{Y}_{n}(x)
       (dg^{(k-1)}(x)-d\tilde{g}^{(k-1)}_{n}(x)) \\
& = & \frac{1}{2}(-1)^{k}
\int _{0}^{\infty} 
(G_{k}(x)-\tilde{H}_{n}(x))(dg^{(k-1)}(x)-d\tilde{g}^{(k-1)}_{n}(x)) 
\\
&&  + \ (-1)^{k} \int _{0}^{\infty}
      (\tilde{H}_{n}(x)-\mathbb{Y}_{n}(x))(dg^{(k-1)}(x)-d\tilde{g}^{(k-1)}_{n}(x)) \\
& \geq &  (-1)^{k}
           \int _{0}^{\infty} (\tilde{H}_{n}(x)-\mathbb{Y}_{n}(x))
        (dg^{(k-1)}(x)-d\tilde{g}^{(k-1)}_{n}(x)).
\end{eqnarray*}
\par \noindent To see that, we notice (using integration by parts) that
\begin{eqnarray*}
(-1)^{k} \int _{0}^{\infty} (G_{k}(x)-\tilde{H}_{n}(x))
          (dg^{(k-1)}(x)-d\tilde{g}^{(k-1)}_{n}(x))
= \int _{0}^{\infty} (g(x)-\tilde{g}_{n}(x))^{2} dx.
\end{eqnarray*}
But condition (\ref{Fenchel}) implies that
\begin{eqnarray*}
\int _{0} ^{\infty} (\tilde{H}_{n}(x)- \mathbb{Y}_{n}(x)) d\tilde{g}_{n}^{(k-1)}(x) =0.
\end{eqnarray*}
Therefore,
\begin{eqnarray*}
Q_{n}(g)- Q_{n}(\tilde{g}_{n})
\geq \int _{0}^{\infty}( \tilde{H}_{n}(x)- 
\mathbb{Y}_{n}(x))(-1)^{k}dg^{(k-1)}(x) \geq 0,
\end{eqnarray*}
since $\tilde{H}_{n} \geq \mathbb{Y}_{n}$ and
$(-1)^{k-2}dg^{(k-1)}(x) = (-1)^{k}dg^{(k-1)}(x) \geq 0 $
because \\ $ (-1)^{k-2}g^{(k-2)}$ is
convex.

Conversely, take $g_{t} \in \mathcal{M}_{k}$ to be
\begin{eqnarray*}
g_{t}(x) = \frac{(t-x)^{k-1}_{+}}{(k-1)!}, \qquad x \geq 0.
\end{eqnarray*}
We have:
\begin{eqnarray*}
\lim_{\epsilon \to 0}
\frac{Q_{n}(\tilde{g}_{n}+\epsilon g_{t}) - Q_{n}(\tilde{g}_{n})}{\epsilon}
& =& \int _{0}^{t}\frac{(t-x)^{k-1}}{(k-1)!} \tilde{g}_{n}(x) dx
       -  \int _{0}^{t} \frac{(t-x)^{k-1}}{(k-1)!} d\mathbb{G}_{n}(x).
\end{eqnarray*}
Using integration by parts, we obtain
\begin{eqnarray*}
0 \le \lim_{\epsilon \to 0}
\frac{Q_{n}(\tilde{g}_{n}+\epsilon g_{t}) - Q_{n}(\tilde{g}_{n})}{\epsilon}
= \tilde{H}_{n}(t)-\mathbb{Y}_{n}(t) \, .
\end{eqnarray*}
Finally, since $\tilde{g}_{n}$ maximizes $Q_{n}$ it follows that
\begin{eqnarray*}
0 & = & \lim_{\epsilon \to 0}
\frac{Q_{n}((1+\epsilon)\tilde{g}_{n}) - Q_{n}(\tilde{g}_{n})}{\epsilon}
  =  \int_0^{\infty} \tilde{g}_n ^2 (x) dx
      - \int_0^{\infty} \tilde{g}_n (x) d \GG_n (x) \\
& = & \int_{0}^{\infty} (\tilde{H}_{n}(x)-\mathbb{Y}_{n}(x))(-1)^{k-1}
        d \tilde{g}_n^{(k-1)}(x),
\end{eqnarray*}
which holds if and only if the equality in (\ref{Fenchel}) holds.
\hfill$\blacksquare$
\medskip

In order to prove that the LSE is a spline of degree $k-1$, we need the following result.

\par \noindent 
\begin{lemma}
\label{Zeros} 
Let $\lbrack a,b \rbrack \subseteq (0,\infty)$ and let $g$ be a nonnegative and nonincreasing function on $\lbrack a, b \rbrack$. For any polynomial $P_{k-1}$ of degree $ \leq k-1$ on $\lbrack a, b \rbrack$, if the function
\begin{eqnarray*}
\Delta(t) = \int_0^t (t-s)^{k-1}g(s) ds - P_{k-1}(s), \ \ t \in \lbrack a,b \rbrack
\end{eqnarray*} 
admits infinitely many zeros in $\lbrack a,b \rbrack$, then there exists $t_0  \in \lbrack a,b \rbrack$ such that $g \equiv 0$ on $\lbrack t_0, b \rbrack$ and $g > 0$ on $\lbrack a, t_0)$ if $ t_0 > a$.  
\end{lemma}
\medskip

\par \noindent 
\textbf{Proof.} 
By applying the mean value theorem $k$ times, it follows that $(k-1)! g = \Delta^{(k)}  $ admits infinitely many zeros in $\lbrack a,b \rbrack$. But since $g$ is assumed to be nonnegative and nonincreasing, this implies that if $t_0$ is the smallest zero of $g$ in $\lbrack a,b \rbrack$, then  $g \equiv 0$ on $\lbrack t_0, b \rbrack$. By definition of $t_0$, $g > 0$ on $\lbrack a, t_0)$ if $t_0 > a$. \hfill $\blacksquare$

\medskip

\par \noindent 
\begin{remark} In the previous lemma, the assumption that $\Delta$ has infinitely many zeros can be weakened. Indeed, we obtain the same conclusion if we assume that $\Delta$ has $k+1$ distinct zeros in $\lbrack a, b \rbrack$.
\end{remark}

Now, we will use the characterization of the LSE \ $\tilde{g}_n$ together with the previous lemma to show that it is a finite mixture of $Beta(1,k)$'s. We know from Proposition \ref{Fenchel} that $\tilde{g}_n$ is the LSE if and only if 
\begin{eqnarray}\label{Ineq}
\tilde{H}_n(t) \geq \mathbb{Y}_n(t), \ \ \textrm{ for $t > 0$},
\end{eqnarray}
and
\begin{eqnarray}\label{Eq}
\int_0^\infty \left(\tilde{H}_n(t) - \mathbb{Y}_n(t) \right) d\tilde{g}^{(k-1)}_n(t) = 0
\end{eqnarray}
where
\begin{eqnarray*}
\tilde{H}_n(t) = \int_0^t \frac{(t-s)^{k-1}}{(k-1)!} \tilde{g}_n(t)dt,
\end{eqnarray*}
and
\begin{eqnarray*}
\mathbb{Y}_n(t) = \int_0^t \frac{(t-s)^{k-1}}{(k-1)!} d\mathbb{G}_n(t).
\end{eqnarray*}
The condition in (\ref{Eq}) implies that $\tilde{H}_n$ and $\mathbb{Y}_n$ have to be equal at any point of increase of the monotone function $(-1)^{k-1} \tilde{g}^{(k-1)}_n$. Therefore, the set of points of increase of $(-1)^{k-1}\tilde{g}^{(k-1)}_n$ is included in the set of zeros of the function $\tilde{\Delta}_n = \tilde{H}_n - \mathbb{Y}_n$. Now, note that $\mathbb{Y}_n$ can be given by the explicit expression:
\begin{eqnarray*}
\mathbb{Y}_n(t) = \frac{1}{(k-1)!} \frac{1}{n} \sum_{j=1}^n (t - X_{(j)})^{k-1}_{+}, \ \ \textrm{ for $t > 0$}.
\end{eqnarray*}
In other words, $\mathbb{Y}_n$ is a spline of degree $k-1$ with simple knots $X_{(1)}, \cdots, X_{(n)}$ (for a definition of the multiplicity 
of knots, see e.g. \mycite{deboor:78}, page 96,
 or \mycite{devore-lorentz:93}, page 140). 
Also note that the function $(-1)^{k-1} \tilde{g}^{(k-1)}_n$ 
cannot have a positive density  
 with respect to Lebesgue measure $\lambda$. Indeed, if we assume otherwise, then we can find $ 0 \leq j \leq n$ and an interval $I \subset (X_{(j)}, X_{(j+1)})$ (with $X_{(0)} = 0$ and $X_{(n+1)} = \infty$) such that $I$ has a nonempty interior, and $\tilde{H}_n \equiv \mathbb{Y}_n $ on $I$. This implies that $\tilde{H}^{(k)}_n \equiv \mathbb{Y}^{(k)}_n \equiv 0$, since $\mathbb{Y}_n$ is a polynomial of degree $k-1$ on $I$, and hence $\tilde{g}_n \equiv 0$ on I. But the latter is impossible since it was assumed that $(-1)^{k-1} \tilde{g}^{(k-1)}_n $ was strictly increasing on $I$. Thus the monotone function $(-1)^{k-1} \tilde{g}^{(k-1)}_n$ can have only two components: discrete and singular. In the following theorem, we will prove that it is actually discrete with finitely many points of jump.

\medskip

\par \noindent 
\begin{prop}
\label{Discretness of the LSE} 
There exists $m \in \mathbb{N}\backslash \{0\}$, $\tilde{a}_1, \cdots, \tilde{a}_m$ and $\tilde{w}_1, \cdots, \tilde{w}_m$  such that for all $x > 0$, the LSE \ $\tilde{g}_n$ is given by
\begin{eqnarray}\label{FormLSE}
\tilde{g}_n(x) = \tilde{w}_1 \frac{k(\tilde{a}_1 -x)^{k-1}_{+}}{\tilde{a}^k_1} + \cdots + \tilde{w}_m \frac{k(\tilde{a}_m -x)^{k-1}_{+}}{\tilde{a}^k_m}. 
\end{eqnarray}
\end{prop}

\medskip

\par\noindent 
\textbf{Proof.} We need to consider two cases: 
\medskip
\par \noindent (i) The number of zeros of $\tilde{\Delta}_n = \tilde{H}_n - \mathbb{Y}_n$ is finite. This implies by (\ref{Eq}) that the number of points of increase of $(-1)^{k-1}\tilde{g}^{(k-1)}_n$ is also finite. Therefore, $(-1)^{k-1} \tilde{g}^{(k-1)}_n$ is discrete with finitely many jumps and hence $\tilde{g}_n$ is of the form given in (\ref{FormLSE}). 
\medskip
\par \noindent 
(ii) Now, suppose that $\tilde{\Delta}_n$ has infinitely many zeros. Let $j$ be the smallest integer in $\{0, \cdots, n-1\}$ such that $\lbrack X_{(j)}, X_{(j+1)} \rbrack$ contains infinitely many zeros of $\tilde{\Delta}_{n}$ (with $X_{(0)} = 0$ and $X_{(n+1)} = \infty$).  By Lemma \ref{Zeros}, if $t_j $ is the smallest zero of $\tilde{g}_n$ in $\lbrack X_{(j)}, X_{(j+1)} \rbrack$, then $\tilde{g}_n \equiv 0$ on $\lbrack t_j, X_{(j+1)} \rbrack$ and $\tilde{g}_n > 0$ on $\lbrack X_{(j)}, t_j)$ if $t_j > X_{(j)}$. Note that from the proof of Proposition \ref{ExistenceLSE}, we know that the minimizing measure $\tilde{\mu}_n$ does not put any mass on $(0, X_{(1)} \rbrack$, and hence the integer $j$ has to be strictly greater than 0. 

Now, by definition of $j$, $\tilde{\Delta}_n$ has finitely many zeros to the left of $X_{(j)}$, 
which implies that  $(-1)^{k-1} \tilde{g}^{(k-1)}_n$ has finitely many 
points of increase in $(0, X_{(j)})$. We also know that 
$\tilde{g}_n \equiv 0$ on $\lbrack t_j, \infty)$. Thus we only need to 
show that the number of points of increase of $(-1)^{k-1} \tilde{g}^{(k-1)}_n$ 
in $\lbrack X_{(j)}, t_j )$ is finite, when $t_j > X_{(j)}$. This can be argued 
as follows: Consider $z_j$ to be the smallest zero of $\tilde{\Delta}_n$ in 
$\lbrack X_{(j)}, X_{(j+1)})$. If $z_j \geq t_j$, then we cannot possibly have 
any point of increase of $(-1)^{k-1} \tilde{g}^{(k-1)}_n$ in $\lbrack X_{(j)}, t_j)$ 
because it would imply that we have a zero of $\tilde{\Delta}_n$ that is strictly 
smaller than $z_j$. If $z_j < t_j$, then for the same reason, $(-1)^{k-1} \tilde{g}^{(k-1)}_n$ 
has no point of increase in $\lbrack X_{(j)}, z_j)$. Finally, 
$(-1)^{k-1} \tilde{g}^{(k-1)}_n$ cannot have infinitely many points of increase in 
$\lbrack z_j, t_j)$ because that would imply that $\tilde{\Delta}_n$ has  infinitely 
zeros in $(z_j, t_j)$, and hence by Lemma \ref{Zeros}, we can find 
$t'_j \in (z_j, t_j)$ such that $\tilde{g}_n \equiv 0$ on $\lbrack t'_j, t_j \rbrack$. 
But this impossible since $\tilde{g}_n > 0$ on $\lbrack X_{(j)}, t_j)$.   
\hfill $\blacksquare$    
\medskip

\begin{remark} 
{\rm We have not succeeded in extending Proposition~\ref{ExistenceLSE}  
to the case $k=\infty$.  It is possible to prove the existence of 
a least squares estimator if the maximization is carried 
over over ${\cal D}_{\infty} \cap L_2(\lambda)$ rather than 
${\cal M}_{\infty} \cap L_2 (\lambda)$, but this does not seem 
(to us) to be the right direction to proceed.}
\end{remark}

\section{Consistency}
\setcounter{equation}{0}

In this section, we will prove that both the MLE and
LSE are strongly consistent. Furthermore,
we will show that this consistency is uniform on
intervals of the form $\lbrack c, \infty)$, where $c >0$.

\subsection{Consistency of the maximum likelihood estimator}

Consistency of the maximum likelihood estimators for the classes 
${\cal D}_k$ in the sense of Hellinger convergence of the mixed density
is a relatively simple straightforward consequence of the methods 
of \mycite{pfanz:88},  \mycite{vdg:93}, and \mycite{vdg:96}. 
As usual, the Hellinger distance 
$H$ is given by $H^2 (p,q) = (1/2) \int \{ \sqrt{p} - \sqrt{q})^2 d \mu$
for any common dominating measure $\mu$.
\medskip

\begin{prop}
\label{HellingerConsistMLE}
Suppose that $\hat{g}_{n,k}$ is the MLE of $g_0$ in the class ${\cal D}_k$, 
$k \in \{1, \ldots , \infty\}$.  Then 
$$
H( \hat{g}_{n,k} , g_0) \rightarrow_{a.s.} 0 \qquad \mbox{as} \ \ n\rightarrow \infty \, .
$$
Furthermore $\hat{F}_{n,k} \rightarrow_d F_0$ almost surely 
where $\hat{F}_{n,k} $ is the MLE
of the mixing distribution function $F_0$.
\end{prop}

\medskip
\par\noindent
{\bf Proof.}
This follows from the methods of \mycite{pfanz:88}, 
\mycite{vdg:93}, and \mycite{vdg:96}, by using the Glivenko-Cantelli 
preservation theorems of \mycite{vdvw:00}.  See also 
\mycite{vdg:99}, page 54, example 4.2.4,
and \mycite{well:03b}, pages 
98 to 99.\\ \hfill $\blacksquare$
\medskip

The following lemma establishes a useful bound for $k$-monotone densities.

\medskip

\par\noindent 
\begin{lemma}\label{Bound} If $g$ is a $k$-monotone density function for $k\ge 2$, then
\begin{eqnarray*}
g(x) \leq \frac{1}{x}\left (1-\frac{1}{k} \right )^{k-1}
\end{eqnarray*}
for all $x > 0$.
\end{lemma}

\medskip

\par \noindent \textbf{Proof.} We have 
\begin{eqnarray*}
g(x)
& = & \int_x^{\infty} \frac{k}{y^k} (y-x)^{k-1} dF(y)
       =  \frac{1}{x} \int_x^{\infty}
          \frac{kx}{y} (1 - \frac{x}{y} )^{k-1}  dF(y) \\
& \le & \frac{1}{x} \sup_{x \le y < \infty } \frac{kx}{y} \left (
           1- \frac{x}{y} \right )^{k-1}
       =  \frac{k}{x} \sup_{0< u \le 1} u (1-u)^{k-1} \\
& = & \frac{1}{x}\left (1-\frac{1}{k} \right )^{k-1} \,
\end{eqnarray*}
since, with $g_k (u) = u(1-u)^{k-1}$ we have
$$
g_k'(u) = (1-u)^{k-1} - u(k-1) (1-u)^{k-2} = (1-u)^{k-2} (1-ku)
$$
which equals zero if $u=1/k$ and this yields a maximum. (Note that when $k=2$, this bound equals $1/(2x)$ which agrees with
the bound given by \mycite{jongb:95}, page 117 in this case.)    \hfill $\blacksquare$

\medskip

\par\noindent \begin{prop}\label{consistencyMLE} Let $g_0$ be a
$k$-monotone density on $(0,\infty)$ and fix $c > 0$. Then
\begin{eqnarray*}
\sup_{x \geq c} \vert \hat{g}_n(x)- g_0(x)
\vert \to_{a.s.} 0, \hspace{0.5cm}
\textrm{as} \hspace{0.3cm} n \to \infty.
\end{eqnarray*}

\end{prop}
\medskip

\par\noindent
\textbf{Proof.}
Let $F_0$ be the mixing distribution function
associated with $g_0$. Then for all $x > 0$, we
have
\begin{eqnarray*}
g_0(x) = \int_0 ^{\infty}  \frac{k (t-x)^{k-1}_{+}}{t^k}dF_0(t).
\end{eqnarray*}

Now, let $Y_1, \cdots, Y_m$ be i.i.d. from $F_0$. Taking $m=n$,
let $\mathbb{F}_n$ be the corresponding empirical distribution and
$g_n$ the mixed density
\begin{eqnarray*}
g_n(x) = \int_0 ^{\infty} \frac{k  (t-x)^{k-1}_{+}}{t^k}
   d\mathbb{F}_n(t), \hspace{0.2cm} x > 0.
\end{eqnarray*}

Let $d > 0 $.  Using integration by parts, we have for all $x > d$
\begin{eqnarray*}
\lefteqn{\vert g_n(x) - g_0(x) \vert} \\
&=& \left \vert \int_x ^{\infty} k
       \frac{(t-x)^{k-1}}{t^k} d(\mathbb{F}_n - F_0)(t) \right \vert \\
& =& \left \vert \int_x ^{\infty}
        k\frac{(k-1)t^k(t-x)^{k-2}-kt^{k-1}(t-x)^{k-1}}{t^{2k}}
       (\mathbb{F}_n-F_0)(t) dt \right \vert \\
&\leq & \left( \int_x ^{\infty} k^2 \frac{(t-x)^{k-2}}{t^{k}}dt
         + \int_x ^{\infty} k^2 x
          \frac{(t-x)^{k-2}}{t^{k+1}}dt \right) \Vert
           \mathbb{F}_n - F_0 \Vert_{\infty} \\
& \leq & \left( \int _d ^{\infty} k\frac{(t-d)^{k-2}}{t^k}dt
           + k^2 \int _d ^{\infty}
           \frac{(t-d)^{k-2}}{t^k}dt \right)
           \Vert \mathbb{F}_n - F_0 \Vert_{\infty} \\
& \leq & \left(
2 k^2 \int _d ^{\infty} \frac{(t-d)^{k-2}}{t^k}dt \right)
         \Vert \mathbb{F}_n - F_0
            \Vert_{\infty} \\
& = & C_d  \Vert \mathbb{F}_n - F_0 \Vert_{\infty}.
\end{eqnarray*}
By the Glivenko-Cantelli theorem, the sequence of
$k$-monotone densities $\left(g_n \right)_n$ satisfies
\begin{eqnarray*}
\sup_{x \in \lbrack d,\infty)}\vert g_n(x) - g_0(x) \vert
\to _{a.s.} 0, \hspace{0.5cm} \textrm{as} \hspace{0.3cm} n \to \infty.
\end{eqnarray*}
 Since the MLE $\hat{g}_n$ maximizes the criterion function over the class $\mathcal{M}_k \cap L_1(\lambda)$, we have
\begin{eqnarray*}
\lim_{\epsilon \searrow 0} \frac{1}{\epsilon}
\left(\psi_n((1-\epsilon)\hat{g}_n
+ \epsilon g_n) - \psi _n(\hat{g}_n) \right) \leq 0,
\end{eqnarray*}
and this is equivalent to
\begin{eqnarray}\label{ineqnMLE}
\int _0 ^\infty \frac{g_n(x)}{\hat{g}_n(x)}d\mathbb{G}_n(x) \leq 1.
\end{eqnarray}
Let $\hat{F}_n$ denote again the MLE of the mixing distribution. By the Helly-Bray theorem, there exists a subsequence $\{\hat{F}_{l}\}$
that converges weakly to some distribution function $\hat{F}$
and hence for all $x > 0$
\begin{eqnarray*}
\hat{g}_{l} (x) \to \hat{g}(x), \ \ \ \textrm{as} \hspace{0.2cm} l \to \infty,
\end{eqnarray*}
where
\begin{eqnarray*}
\hat{g}(x)
= \int_0^{\infty} k
     \frac{(t-x)^{k-1}_{+}}{t^k}d\hat{F}(t), \qquad x > 0.
\end{eqnarray*}
The previous convergence is uniform on intervals of the form
$\lbrack d, \infty)$, $d > 0$. This follows since
$\hat{g}_l$ and $\hat{g}$ are monotone and
$\hat{g}$ is continuous.

Much of the following is along the lines of
\mycite{jongb:95}, pages 117-119, and
\mycite{gjw:01b}, pages 1674-1675.
We are going to show that $\hat{g}$ and the true
density $g_0$ have to be the same.
For $0< \alpha < 1$ define $\eta_{\alpha} = G_0^{-1} (1-\alpha)$.
Fix  $\epsilon $ so small that $\epsilon < \eta_{\epsilon}$.
By (\ref{ineqnMLE})
there is a number $D_\epsilon > 0$ such that
$\hat{g}_l(\eta_{\epsilon}) \geq D_{\epsilon}$
for sufficiently large $l$.
To see this, note that (\ref{ineqnMLE}) implies that
\begin{eqnarray*}
1  \ge  \int_0 ^{\infty} \frac{g_l(x)}{\hat{g}_l(x)}d\mathbb{G}_l(x)
  \ge  \int_{\eta_{\epsilon}}^{\infty}
           \frac{g_l(x)}{\hat{g}_l(x)}d\mathbb{G}_l(x)
  \ge  \frac{1}{ \hat{g}_l (\eta_{\epsilon})} \int_{\eta_{\epsilon} }^{\infty}
            g_l (x) d \GG_l (x) \, ,
\end{eqnarray*}
and hence
\begin{eqnarray*}
\liminf_l  \hat{g}_l (\eta_{\epsilon} ) \ge
\liminf_l \int_{\eta_{\epsilon}}^{\infty} g_l (x) d \GG_l (x)
= \int_{\eta_{\epsilon}}^{\infty} g_0 (x) d G_0 (x) > 0 \, ,
\end{eqnarray*}
by the choice of $\eta_{\epsilon}$
and hence we can certainly take
$D_{\epsilon} = \int_{\eta_{\epsilon}}^{\infty} g_0 (x) d G_0 (x)/2$.

Hence, by continuity of $g_l$ and the bound in Lemma \ref{Bound}
$$
\hat{g}_l (z) \le \frac{1}{z} (1 - \frac{1}{k} )^{k-1} \equiv 
\frac{e_k}{z} \, ,\qquad
g_l (z) \le \frac{1}{z} (1 - \frac{1}{k} )^{k-1}\equiv \frac{e_k}{z} \, ,
$$
$g_l/\hat{g}_l$ is uniformly bounded on the
interval $\lbrack \epsilon, \eta_\epsilon \rbrack $.
That is, there exist two constants
$\underline{c}_\epsilon$ and $\overline{c}_\epsilon$ such that for all
$x \in \lbrack \epsilon, \eta_\epsilon \rbrack $
\begin{eqnarray*}
\underline{c}_\epsilon
\leq \frac{g_l(x)}{\hat{g}_l(x)} \leq \overline{c}_\epsilon.
\end{eqnarray*}
In fact,
\begin{eqnarray*}
  \frac{g_l(x)}{\hat{g}_l(x)} \le \frac{g_l (\epsilon)}{\hat{g}_l 
(\eta_{\epsilon})}
\le \frac{\epsilon^{-1} e_k}{D_\epsilon},
\end{eqnarray*}
while
\begin{eqnarray*}
\frac{g_l(x)}{\hat{g}_l(x)} \ge \frac{g_l (\eta_{\epsilon} 
)}{\hat{g}_l (\epsilon)}
\ge \frac{g_0 (\eta_\epsilon)/2}{\epsilon^{-1} e_k}
\end{eqnarray*}
using the (uniform) convergence of $g_l $ to $g_0$. Therefore
\begin{eqnarray*}
\frac{g_l(x)}{\hat{g}_l(x)} \to \frac{g_0(x)}{\hat{g}(x)}
\end{eqnarray*}
uniformly on $\lbrack \epsilon, \eta_\epsilon]$. For sufficiently large $l$,
we have using (\ref{ineqnMLE})
\begin{eqnarray*}
\int_\epsilon ^{\eta_\epsilon}
    \frac{g_0(x)}{\hat{g}(x)}d \GG_l(x) \leq \int_\epsilon 
^{\eta_\epsilon}\left(
        \frac{g_l(x)}{\hat{g}_l(x)}
       + \epsilon \right) d \GG_l(x) \leq 1+ \epsilon.
\end{eqnarray*}
But since $\GG_l$ converges weakly to $G_0$ the
distribution function of $g_0$ and $g_0/\hat{g}$
is continuous and bounded on $\lbrack
\epsilon, \eta_\epsilon \rbrack$, we conclude that
\begin{eqnarray*}
\int_\epsilon ^{\eta_\epsilon}\frac{g_0(x)}{\hat{g}(x)}dG_0(x)
\leq 1+\epsilon.
\end{eqnarray*}
Now, by Lebesgue's monotone convergence theorem, we conclude that
\begin{eqnarray*}
\int_0 ^{\infty}\frac{g_0(x)}{\hat{g}(x)}dG_0(x) \leq 1,
\end{eqnarray*}
which is equivalent to
\begin{eqnarray}\label{Resineqn}
\int_0 ^{\infty}\frac{g^2_0(x)}{\hat{g}(x)}dx \leq 1.
\end{eqnarray}
Define $\tau = \int_0^\infty \hat{g}(x)dx$.
Then $\hat{h} = \tau^{-1}\hat{g}$ is a $k$-monotone density.
By (\ref{Resineqn}), we have that
\begin{eqnarray*}
\int_0^{\infty} \frac{g^2_0(x)}{\hat{h}(x)}dx
= \tau \int_0 ^{\infty} \frac{g^2_0(x)}{\hat{g}(x)}dx \leq \tau.
\end{eqnarray*}
Now consider the function
\begin{eqnarray*}
K(g)= \int_0^{\infty} \frac{g^2_0(x)}{g(x)}dx
\end{eqnarray*}
defined on the class $\mathcal{C}_d$ of all continuous densities
$g$ on $\lbrack 0, \infty)$. Minimizing $K$ is equivalent to minimizing
\begin{eqnarray*}
\int_0^{\infty} \left(\frac{g^2_0(x)}{g(x)} + g(x)\right)dx.
\end{eqnarray*}
It is easy to see that the integrand is
minimized pointwise
by taking $g(x)=g_0(x)$. Hence $\inf_{\mathcal{C}_d} K(g) \geq 1$. In
particular,
$K(\hat{h}) \geq 1$ which implies that $\tau=1$.
Now, if $g \ne g_0$ at a point $x$, it follows
that $g \ne g_0$ on an interval of positive length.
Hence, $g_0 \ne g \Rightarrow K(g) > 1$.
We conclude that we have necessarily $\hat{h}= \hat{g} = g_0$.

We have proved that from each subsequence of
$\hat{g}_n$, we can extract a further subsequence that converges
to $g_0$ almost surely. The
convergence is again uniform on intervals of the form
$\lbrack c, \infty)$, $c > 0$ by
monotonicity of $\hat{g}_n$ and $\hat{g}$ and continuity of $g_0$.
\hfill $\blacksquare$
\medskip

\par\noindent 
\begin{corollary}
\label{uniformconsistDerivMLE} 
Let $c > 0$.
For $j=1, \cdots, k-2$,
\begin{eqnarray*}
\sup_{x \in \lbrack c, \infty)} \vert \hat{g}^{(j)}_n(x) - g^{(j)}_0(x)
\vert \to_{a.s.} 0, \hspace{0.2cm} \textrm{as} \hspace{0.2cm} n \to \infty,
\end{eqnarray*}
and for each $x > 0$ at which $g_0$ is $k-1$-times differentiable,
\begin{eqnarray*}
  \hat{g}^{(k-1)}_n(x) \rightarrow_{a.s.} g^{(k-1)}_0(x) \, .
\end{eqnarray*}
\end{corollary}
\medskip

\par\noindent
\textbf{Proof.}
This follows along the lines of the proof
in  \mycite{jongb:95}, page 119, and
\mycite{gjw:01b}, Lemma 3.1, page 1675.
\hfill $\blacksquare$
\medskip

\subsection{The Least Squares estimator}
\medskip

We also have strong and uniform consistency of the LSE
$\tilde{g}$ on intervals of the form
$ \lbrack c,\infty), c > 0$.
\medskip

\par\noindent \begin{prop}\label{consistencyLSE}
Fix $c > 0$ and suppose that the true $k$-monotone density  $g_0$ satisfies
$\int_0^{\infty} x^{-1/2} d G_0 (x) < \infty$. Then \ $\| \tilde{g}_n - g_0 \|_2 \rightarrow_{a.s.} 0$, and  
\begin{eqnarray*}
\sup_{x \geq c} \vert \tilde{g}_n(x)- g_0(x) \vert
\to_{a.s.} 0, \hspace{0.2cm} \textrm{as} \hspace{0.2cm} n \to \infty.
\end{eqnarray*}
\medskip
\end{prop}
\par\noindent \textbf{Proof.}
The main difficulty here is that we don't know whether the
LSE $\tilde{g}_n$ is a genuine density; i.e. $\tilde{g}_n \in {\cal M}_k$
but not necessarily $\tilde{g}_n \in {\cal D}_k$.
Once we show  that
$\tilde{g}_n$ stays bounded in $L_2$ with high probability,
the proof of consistency will be much like the one used for $k=2$; 
i.e., consistency of the
LSE of a convex and decreasing density (see \mycite{gjw:01b}). The
proof for $k=2$ is based on the very important fact that the LSE is a 
density, which helps in
showing that $\tilde{g}_n$ at the last jump point
$\tau_{n} \in \lbrack 0, \delta \rbrack $  of $\tilde{g}'_n$
for a fixed $\delta > 0$ is uniformly bounded. The proof would have
been similar if we only  knew that
\begin{eqnarray*}
\int_{0}^{\infty}\tilde{g}_n(x) dx = O_p (1) \, .
\end{eqnarray*}

Here we will first show that
$\int_0^{\infty} \tilde{g}_n^2 d\lambda = O(1) $ almost surely.
 From the last display in the proof of Proposition \ref{characterizationLSE}
\begin{eqnarray*}
\int_0^{\infty}\tilde{g}_n^2(x)dx
= \int_0^{\infty}\tilde{g}_n(x)d\mathbb{G}_n(x)
\end{eqnarray*}
and hence
\begin{eqnarray}
\sqrt{\int_0^{\infty}\tilde{g}_n^2(x)dx}
= \int_0^{\infty}\tilde{u}_n(x)d\mathbb{G}_n(x),
\label{L2NormEquality}
\end{eqnarray}
where $\tilde{u}_n \equiv \tilde{g}_n /\| \tilde{g}_n\|_2$
satisfies $\Vert \tilde{u}_n \Vert _2  =1$.
Take $\mathcal{F}_k $ to be the class of functions
\begin{eqnarray*}
\mathcal{F}_k = \bigg  \{ g \in \mathcal{M}_k,
   \int_0^\infty g^2 d \lambda =1  \bigg \}.
\end{eqnarray*}
In the following, we show that
$\mathcal{F}_k$ has an envelope $G \in L_1(G_0)$. \\

Note that for $g \in {\cal F}_k$ we have
\begin{eqnarray*}
1 = \int_0^{\infty} g^2 d \lambda \ge \int_0^x g^2 d \lambda
\ge x g^2 (x) \, ,
\end{eqnarray*}
since $g$ is decreasing.  Therefore
\begin{eqnarray*}
g(x) \le \frac{1}{\sqrt{x}}  \equiv G(x)
\end{eqnarray*}
for all $x > 0$ and $g \in {\cal F}_k$;
i.e. $G$ is an envelope for the class $\mathcal{F}_k$.
Since $G \in L_1 (G_0)$ (by our hypothesis)
it follows from the strong law that
\begin{eqnarray*}
\int_0^\infty \tilde{u}_n(x) d\mathbb{G}_n(x)
      \leq \int_0^\infty G(x)d\mathbb{G}_n(x)
        \to_{a.s.} \int_0^\infty
G(x)d{G}_0(x), \hspace{0.3cm} \textrm{as} \hspace{0.3cm} n \to \infty
\end{eqnarray*}
and hence by (\ref{L2NormEquality}) the integral
$\int_0^\infty \tilde{g}^2_n d\lambda $ is
bounded (almost surely)  by some constant $M_k$.

Now we are ready to complete the proof.
Most of the following arguments are similar to those of proof
of consistency of the LSE when $k=2$ as given in \mycite{gjw:01b}.
\medskip

Let $\delta > 0$ and $\tau_{n}$ be the last jump
point of $\tilde{g}^{(k-1)}_n$ if there are jump points
in the interval $(0,\delta \rbrack$,
otherwise we take $\tau_{n}$ to be 0. To show that the sequence
$\left(\tilde{g}_n(\tau_n)\right)_n$ stays bounded, we consider two cases:
\begin{enumerate}
\item $\tau_n \geq \delta /2$. Let $n$ be large enough so that
$\int_0^{\infty} \tilde{g}_n^2 d \lambda \le M_k$.
We have
\begin{eqnarray}\label{Bounddelta}
\tilde{g}_n(\tau_n) & \leq & \tilde{g}_n(\delta /2)
  \leq  (2/\delta)(\delta/2) \tilde{g}_n(\delta /2)
  \leq  (2/\delta)\int_0^{\delta/2}\tilde{g}_n(x)dx \nonumber \\
&\leq & (2/\delta) \sqrt{\delta/2}\sqrt{\int_0^{\delta /2}\tilde{g}_n^2(x)dx}
  \leq  \sqrt{2/\delta}\sqrt{\int_0^{\infty}\tilde{g}_n^2(x)dx} \nonumber\\
& = &  \sqrt{2M_k/\delta}.
\end{eqnarray}
\item $\tau_n < \delta /2$.  We have
\begin{eqnarray*}
\int_{\tau_n}^{\delta}\tilde{g}_n(x)dx & \leq & \sqrt{\delta -\tau_n} 
\sqrt{\int_{\tau_n}^{\delta}\tilde{g}_n^2(x)dx}\\
& \leq & \sqrt{\delta}\sqrt{\int_{0}^{\infty}\tilde{g}_n^2(x)dx}
  =   \sqrt{\delta M_k}.
\end{eqnarray*}
\par \noindent 
Using the fact that $\tilde{g}_n$ is a polynomial
of degree $k-1$ on the interval $\lbrack \tau_{n}, \delta \rbrack $ we have
\begin{eqnarray*}
\lefteqn{
\sqrt{\delta M_k}
\geq  \int_{\tau_{n}}^{\delta}\tilde{g}_n(x)dx }\\
& = & \tilde{g}_n(\delta)(\delta - \tau_{n})-
       \frac{\tilde{g}'_n(\delta)}{2}(\delta-\tau_{n})^2 
 +  \cdots +(-1)^{k-1}\frac{\tilde{g}^{(k-1)}_n(\delta)}{k!}(\delta
       - \tau_{n})^k \\
&\geq  &  (\delta - \tau_{n})\left(\tilde{g}_n(\delta)
           +\frac{1}{k} 
           (-1) \tilde{g}'_n(\delta)(\delta-\tau_{n})
       \right . \\
&& \qquad \qquad  \left . 
+ \ \cdots+(-1)^{k-1}
           \frac{\tilde{g}^{(k-1)}_n(\delta)}{(k-1)!}(\delta
           - \tau_{n})^{k-1} \right) \\
& = & (\delta - \tau_{n})\left(\tilde{g}_n(\delta)\left(1
        -\frac{1}{k} \right)
        + \frac{1}{k}\tilde{g}_n(\tau_{n}) \right)\\
&\geq & \frac{\delta}{2k} \tilde{g}_n(\tau_{n})
\end{eqnarray*}
\end{enumerate}
and hence
$\tilde{g}_n(\tau_{n}) \leq 2k \sqrt{M_k / \delta}$.
Therefore, combining the bounds, we have for large $n$
\begin{eqnarray}\label{BoundLSE}
\tilde{g}_n(\tau_n) \leq 2k \sqrt{M_k / \delta} = C_k.
\end{eqnarray}
\medskip
Now, since $\tilde{g}_{n}(\delta) \leq \tilde{g}_{n}(\tau_{n})$,
the sequence $\tilde{g}_{n}(x)$  is uniformly bounded
almost surely for all $x \geq \delta$.
Using a Cantor diagonalization argument, we can find a
subsequence $\{ n_l\} $ so that, for
each $x \geq \delta$, $g_{n_l}(x) \to \tilde{g}(x)$, as $l \to \infty$.
By Fatou's lemma, we
have
\begin{eqnarray}\label{EssIneq}
\int_{\delta}^\infty (\tilde{g}(x) - g_0(x))^2 dx
\leq \liminf_{l \to \infty} \int_{\delta}^\infty (\tilde{g}_{n_l} (x) 
- g_0(x))^2 dx.
\end{eqnarray}
On the other hand,  the 
characterization of $\tilde{g}_{n}$ implies that $Q_n (\tilde{g}_n) \le Q_n (g_0)$,
and this yields 
\begin{eqnarray*}
\int_0^\infty (\tilde{g}_{n}(x) - g_0(x))^2 dx \le 2 \int_0^{\infty} (\tilde{g}_n (x) - g_0 (x) ) d ( \GG_n (x) - G_0 (x)) \, .
\end{eqnarray*}
Thus we can write
\begin{eqnarray}
\lefteqn{
\int_\delta^\infty (\tilde{g}_{n_l}(x) - g_0(x))^2 dx
 \leq   \int_0^\infty (\tilde{g}_{n_l}(x) - g_0(x))^2 dx } \nonumber \\
& \leq & 2  \int_0^\infty
           (\tilde{g}_{n_l}(x) - g_0(x))d(\mathbb{G}_{n_l}(x) - G_0(x))
           \to_{a.s.}  0,   \label{Limit}
\end{eqnarray}
as $l \to \infty$.
The last convergence is justified as follows: since $\int_0^\infty
\tilde{g}^2_{n_l}d\lambda$
is  bounded almost surely,  we can find a
constant $C > 0$ such that
$\tilde{g}_{n_l} - g_0$ admits $G(x) = C/\sqrt{x}, x > 0$, as an envelope.
Since $G \in L_1(G_0)$ by hypothesis and since the class of
functions $\{ (g - g_0)1_{[G \le M]} : \ g \in {\cal M}_k \cap L_2 
(\lambda) \}$
is a Glivenko-Cantelli class for every $M>0$ (each element is a difference
of two bounded monotone functions) (\ref{Limit})
holds.
 From (\ref{EssIneq}), we conclude  that
\begin{eqnarray*}
\int_{\delta}^\infty (\tilde{g}(x) - g_0(x))^2 dx \leq 0 \, ,
\end{eqnarray*}
and therefore, $\tilde{g} \equiv g_0$ on $(0,\infty)$ since $\delta > 0$
can be chosen arbitrarily small.
We have proved that there exists $\Omega_0$ with
$P(\Omega_0)=1$ and such that for each $\omega \in \Omega_0$ and any 
given subsequence
$\tilde{g}_{n_k}(\cdot , \omega)$, we can extract a further subsequence
$\tilde{g}_{n_l}(\cdot , \omega)$ that converges to $g_0$ on $(0, \infty)$.
It follows that $\tilde{g}_n$ converges to
$g_0$ on $(0,\infty)$, and this convergence is
uniform on intervals of the form $\lbrack c, \infty)$, $c > 0$
by the monotonicity and continuity of $g_0$.
\hfill $\blacksquare$

\medskip

\par \noindent 
\begin{corollary}
\label{uniformconsistDerivLSE} 
Let $c > 0$. Under the assumption of Proposition \ref{consistencyLSE}, we have 
for $j=1, \cdots, k-2$,
\begin{eqnarray*}
\sup_{x \in \lbrack c, \infty)} \vert \tilde{g}^{(j)}_n(x) - g^{(j)}_0(x)
\vert \to_{a.s.} 0, \hspace{0.2cm} \textrm{as} \hspace{0.2cm} n \to \infty,
\end{eqnarray*}
and for each $x > 0$ at which $g_0$ is $k-1$-times differentiable,
\begin{eqnarray*}
  \tilde{g}^{(k-1)}_n(x) \rightarrow_{a.s.} g^{(k-1)}_0(x) \, .
\end{eqnarray*}
\end{corollary}
\medskip

\par \noindent \textbf{Proof.} See the proof of Corollary \ref{uniformconsistDerivMLE}.   \hfill $\blacksquare$
\bigskip
 
\section{Asymptotic Minimax risk lower bounds for the rates of convergence}
\setcounter{equation}{0}

In this section our goal is to derive minimax lower bounds for the behavior 
of {\sl any estimator} of a $k-$monotone density $g$ and its first $k-1$
derivatives at a point $x_0$ for which the $k-$th derivative exists and
is non-zero.  The proof will rely upon the basic Lemma 4.1 of 
\mycite{gro:96};
see also \mycite{jongb:00}.  This basic method seems to go back to 
\mycite{donoho-liu:87} and \mycite{donoho-liu:91}).   The relationship 
of our results to other
rate results due to \mycite{kiefer:82}, \mycite{stone:80}, \mycite{fan:91}, and 
\mycite{zhang:90} will be discussed later in the section.

As before, let ${\cal D}_k$ denote the class of $k-$monotone densities
on $[0,\infty)$.  Here is the notation we will need.
Consider estimation of the $j-$th derivative of $g \in {\cal D}_k$ at
$x_0$ for $j\in \{ 0, 1 , \ldots , k-1\}$.
If $\hat{T}_n$ is an arbitrary estimator of the
real-valued functional  $T$ of $g$,
then the $(L_1-)$minimax risk based on a sample $X_1 ,\ldots , X_n$ of size $n$
from $g$ which is known to be in a suitable subset ${\cal D}_{k,n}$ 
of ${\cal D}_k$
is defined by
\begin{eqnarray*}
MMR_1 (n,T,\mathcal{D}_{k,n})
= \inf_{t_n}
\sup_{g \in \mathcal{D}_{k,n}} E_{g} \big \vert \hat{T}_{n}-T g \big \vert \, .
\end{eqnarray*}
Here the infimum ranges over all possible measurable functions
$t_n : \RR^n \rightarrow \RR$, and $\hat{T}_n = t_n (X_1 , \ldots , X_n )$.
When the subclasses ${\cal D}_{k,n}$ are taken to be shrinking to
one fixed $g_0 \in {\cal D}_k$, the minimax risk is called {\sl local}
at $g_0$.
The shrinking classes (parametrized by $\tau>0$) used here are Hellinger
balls centered at $g_0$:
\begin{eqnarray*}
{\cal D}_{k,n} \equiv {\cal D}_{k,n,\tau} =
\bigg \{ g \in \mathcal{D}_{k} :
\  H^2 (g, g_0) = \frac{1}{2} \int_0^{\infty}
( \sqrt{ g (x)} - \sqrt{g_0 (x)} )^2 dx \le \tau/n   \bigg \},
\end{eqnarray*}
The behavior, for $n\rightarrow \infty$ of such a local minimax risk
$MMR_1$ will depend on $n$ (rate of convergence to zero) and the density $g_0$
toward which the subclasses shrink.  The following lemma is the basic tool for
proving such a lower bound.
\medskip

\par\noindent \begin{lemma} Assume that there exists some subset $\{ g_{\epsilon} : \epsilon > 0 \}$
of densities in ${\cal D}_{k,n}$ such that, as $\epsilon \downarrow 0$,
$$
H^2 (g_{\epsilon}, g_0 )\le \epsilon (1+ o(1)) \ \ \mbox{and} \ \
|Tg_{\epsilon} - Tg_0 | \ge (c\epsilon)^r (1+ o (1))
$$
for some $c>0$ and $r>0$.  Then
\begin{eqnarray*}
\sup_{\tau>0} \liminf_{n\rightarrow \infty} \, n^r MMR_1 (n,T, {\cal D}_{k,n} )
\ge \frac{1}{4} \left ( \frac{ cr}{2 e} \right )^r \, .
\end{eqnarray*}

\end{lemma}
\medskip

\par\noindent {\bf Proof.} See \mycite{jongb:95} and \mycite{jongb:00}.
\hfill$\blacksquare$
\medskip

Here is the main result of this section:
\medskip

\par\noindent 
\begin{prop}
\label{MinimaxLowerBoundsForEstimationofDerivativesOfg}
Let $g_0 \in \mathcal{D}_{k}$ and $x_{0}$ be a fixed point in
$(0,\infty)$ such that $g_0$ is $k$ times differentiable
at $x_{0}$ ($k \geq 2$).
An asymptotic
lower bound for the local minimax risk
of any estimator $\hat{T}_{n,j}$ for
estimating the functional
$T_{j}g_0 = g^{(j)}_0(x_{0})$,
is given by:
\begin{eqnarray*}
\sup_{\tau > 0} \liminf_{n \to \infty}
n^{\frac{k-j}{2k+1}}MMR_{1}(n,T_{j},\mathcal{D}_{k,n,\tau}) 
\geq \bigg \{\vert
g^{(k)}_0(x_{0}) \vert ^{2j+1} g_0(x _{0})^{k-j}\bigg \}^{1/(2k+1)}d_{k,j},
\end{eqnarray*}
where $d_{k,j} > 0$, $j \in \{ 0, \ldots ,  k-1\}$.
Here
\begin{eqnarray*}
  d_{k,j} =
\frac{1}{4}
       \bigg(4 \frac{k-j}{2k+1}e^{-1} \bigg)^{\frac{k-j}{2k+1}}
       \frac{\lambda^{(j)}_{k,1}}{\left(\lambda_{k,2}\right)^{\frac{k-j}{2k+1}}}
\end{eqnarray*}
where
\begin{eqnarray*}
\lambda_{k,2} =
2^{4(k+1)} \frac{(2k+3)(k+2)}{(k+1)^2}
\frac{\left((2(k+1))! \right)^2}{ (4k+7)! ((k-1)!)^2
\left( {k\choose k/2-1} \right )^2} ,
\hspace{0.3cm} \textrm{when $k$ is even} \\
\end{eqnarray*}
and
\begin{eqnarray*}
\lambda_{k,2} = 2^{4(k+2)} (2k+3)(k+2)
\frac{\left((2(k+1))! \right)^2}{ (4k+7)!(k!)^2
\left( {k+1 \choose (k-1)/2} \right )^2 }
\hspace{0.3cm} \textrm{when $k$ is odd}
\end{eqnarray*}
and, with $r(x) \equiv (1-x^2)^{k+1} (1+x)$ for $-1 \le x \le 1$ and
$ C_{k,j} \equiv r^{(j)} (0)$,
\begin{eqnarray*}
\lambda^{(j)}_{k,1}
= \left \vert \frac{C_{k,j}}{C_{k,k}} \right \vert, \qquad 0 \leq j \leq k-1 .
\end{eqnarray*}
\end{prop}
\medskip

Proposition ~\ref{MinimaxLowerBoundsForEstimationofDerivativesOfg}
also yields lower bounds for estimation of the 
corresponding mixing distribution function $F$ at a fixed point.

\par\noindent
\begin{corollary}
\label{LowerBoundForMixingDF}
Let $g_0 \in {\cal D}_k$ and let $x_0$ be a fixed point in $(0,\infty)$ such that
$g_0$ is $k-$times differentiable at $x_0$, $k \ge 2$.  Then, for estimating 
$T g_0  = F(x_0)$ where $F_0$ is given in terms of $g_0$ by (\ref{Inverseformula}),
\begin{eqnarray*}
\lefteqn{\sup_{\tau > 0} \liminf_{n \to \infty}
n^{\frac{1}{2k+1}}MMR_{1}(n,T,\mathcal{D}_{k,n,\tau}) }\\
&& \qquad \geq \bigg \{\vert g^{(k)}_0(x_{0}) \vert ^{2k-1} g_0(x _{0})^{}\bigg \}^{1/(2k+1)}
 \frac{x_0^k}{k!}  d_{k,k-1},
\end{eqnarray*}
\end{corollary}

The lower bound results in Proposition~\ref{MinimaxLowerBoundsForEstimationofDerivativesOfg} are 
consistent with the results of \mycite{kiefer:82} and \mycite{stone:80} (although
our result involves a slightly stronger lower bound since the supremum is
over just a local neighborhood of the truth).
In particular, Kiefer showed that rates of 
convergence in estimation cannot be improved by order restrictions, but that 
order restrictions might result in improvements of the constants.  
This latter suggestion has been investigated in  detail in the case of
monotone densities by \mycite{birge:87}, \mycite{birge:89}.  The dependence
of our lower bound on the constants $g_0 (x_0)$ and $g_0^{(k)}(x_0)$
 matches with the known results
for $k=1$ and $k=2$ due to \mycite{gro:85} and \mycite{gjw:01b}, and will 
reappear in the limit distribution theory for $k\ge 3$  in \mycite{balabwell:04c}.

The result of Corollary~\ref{LowerBoundForMixingDF} is consistent with 
the lower bound results of \mycite{zhang:90} and \mycite{fan:91} in the 
deconvolution setting as we now explain.  

To link up with the deconvolution literature 
we transform our scale mixture problem to a location mixture 
or deconvolution problem.  To do this we will reparametrize our $k-$monotone
densities so that the beta kernels converge to the limiting exponential kernels:
Note that if 
$$
g(x) = \int_0^{\infty} \frac{1}{y} \left ( 1 - \frac{y}{kz} \right )_+^{k-1} dF(y) \, ,
$$
then for $X \sim g$, $Z = Z_k \sim k \times \mbox{Beta}(1,k)$, and 
$Y \sim F$ with $Y$ and $Z$ independent, we have 
$$
X \stackrel{d}{=} Z Y\, .
$$
Thus 
$$
X^* \equiv \log X = \log Y + \log Z \equiv Y^* + Z^* \, .
$$
Hence the density $g^*$ of $X^*$ is given 
by 
$$
g^* (x) = \int_{-\infty}^{\infty} \left (1- \frac{1}{k} e^{x-y} \right )_+^{k-1} e^{x-y} dF^* (y)
= \int_{-\infty}^{\infty} f_{Z^*} (x-y) dF^* (y) 
$$
where $F^* (y) = F(e^{y})$ is the distribution function of $Y^*$.

For the completely monotone 
case corresponding to $k=\infty$, the corresponding formulas for $g$ and $g^{*}$ are 
given by
$$
g(x) = \int_0^{\infty} \frac{1}{y} \exp(-x/y) dF(y) \, ,
$$
and 
$$
g^* (x) = \int_{-\infty}^{\infty} \exp(-e^{x-y}) e^{x-y} dF^* (y) = 
\int_{-\infty}^{\infty} f_{Z_{\infty}^*} (x-y) dF^* (y) \, .
$$
\medskip

According to Fan (1991), we need to compute the characteristic 
function $\phi_{Z^*}$ and bound its modulus
 above and below for large arguments.
Thus  we calculate first for $Z_{\infty}^*$:  
from \mycite{abram-steg:64}, page 930, 
\begin{eqnarray*}
\phi_{Z_{\infty}^*} (t) 
 =   \int_{-\infty}^{\infty} e^{itz} e^{-e^z} e^z dz
 =  \int_0^{\infty} e^{it\log v} e^{-v} dv = \Gamma (1+it) \, .
\end{eqnarray*}
Thus by \mycite{abram-steg:64}, page 256, 
$$
| \phi_{Z_{\infty}^*} (t) |^2  = \Gamma(1+it)\Gamma (1-it) 
= \frac{\pi t}{\sinh(\pi t)}
= \frac{2\pi t}{e^{\pi t} - e^{-\pi t}} \, ,
$$
and it follows that 
$$
\sqrt{2 \pi |t|} \exp(- \pi |t|/2) \le 
| \phi_{Y_{\infty}^*} (t) | \le \sqrt{3 \pi |t|} \exp(- \pi |t|/2)
$$
for $|t| \ge 1$.  Thus the hypothesis (1.3) of \mycite{fan:91} holds
with $\beta =1$, $\beta_1 = 1/2$ and $\beta_0 = 1/2$.   This 
implies the first hypothesis of Fan's theorem 4, page 1263, and thus
we are in the case of a ``super-smooth'' convolution kernel.
Fan's second hypothesis is easily satisfied by the current extreme
value distribution function since $f_{Z_\infty^*} (y) = O(|y|^{-2})$
as $y \rightarrow \pm \infty$. 
It therefore follows in the completely monotone case ($k= \infty$) that 
for estimation of $F_0^*(y_0)= F(e^{y_0})$ the resulting minimax lower bound yields 
the rate of convergence $(\log n)^{-1}$.   This rate could also be deduced 
from \mycite{zhang:90}, Corollary 3, page 824. (Note that the tail behavior 
of the characteristic function of our extreme value kernel coincides with the 
tail behavior of the characteristic function of the 
Cauchy kernel and that Zhang's example 2 yields the 
rate $(\log n)^{-1}$ in the case of the Cauchy kernel.)

We can also follow the deconvolution 
approach to obtain a minimax lower bound 
for estimation of the mixing distribution in 
the $k-$monotone case: the characteristic function of 
$Z_k^* = \log Z_k$ is given by 
\begin{eqnarray*}
\phi_{Z_{k}^*} (t) 
& = &  \int_{-\infty}^{\infty} e^{itz} 
           \left (1 - \frac{1}{k}e^{z} \right )_{+}^{k-1} e^z dz
 =  \int_0^{k} e^{it\log v} (1 - v/k)_+^{k-1} dv \\
& = &  \frac{k^{it} \Gamma(k+1) \Gamma (1+it)}{\Gamma (k+1+it)} \, .
\end{eqnarray*}
Thus 
\begin{eqnarray*}
| \phi_{Z_k^*} (t) |^2 
& = &  \frac{k^{it} \Gamma(k+1) \Gamma (1+it)}{\Gamma (k+1+it)} 
             \frac{k^{-it} \Gamma(k+1) \Gamma (1-it)}{\Gamma (k+1-it)}\\
& = & \frac{\Gamma(k+1)^2}{(k+it)(k-1+it) \cdots (1+it) (k-it)(k-1-it)\cdots (1-it)}\\
& = & \frac{(k!)^2}{(k^2 +t^2) \cdots (1+t^2)}   
 \sim  \frac{(k!)^2}{t^{2k}} \qquad \mbox{as} \ \ t \rightarrow \infty\, .
\end{eqnarray*}
It should also be noted  that 
$$
\lim_{k\rightarrow \infty} | \phi_{Z_k^*} (t)|^2
=\lim_{k\rightarrow \infty}  \frac{(k!)^2}{(k^2 +t^2) \cdots (1+t^2)}  = \frac{\pi t}{\sinh(\pi t)} 
=| \phi_{Y_{\infty}^*} (t) |^2 \, .
$$
Thus 
$$
|\phi_{Z_k^*} (t) | \sim \frac{k!}{t^k} \qquad \mbox{as} \ \ t \rightarrow \infty \, ,
$$
and we are in the  situation of a smooth convolution kernel 
of hypothesis (1.4) of Fan (1991), page 1263, 
with Fan's $\beta = k$ in our setting.   Thus Fan's theorem (extended to negative 
values of $l$) gives
our rate of 
convergence for estimating $F^*(y_0)= F(e^{y_0})$ or 
$g^{(k-1)}$  by taking $l=-1$, $\alpha + m = 0$, and $\beta = k$.  
By ``extending'' Fan's theorem further and taking $l=-(k-j)$, we get the 
rate of convergence $n^{-(k-j)/(2k+1)}$, $j=1, \ldots , k-1$ for estimation 
of $g_0^{(j)} (x_0)$.

\medskip

\par\noindent
\textbf{Proof of Proposition~\ref{MinimaxLowerBoundsForEstimationofDerivativesOfg}}. \ \
Let $\mu$ be a positive number and  consider the function
$g_{\mu}$ defined by:
\begin{eqnarray*}
g_\mu(x) = g_0(x) + s(\mu)(x_0+\mu -x)^{k+1}(x-x_0 + \mu)^{k+2}
1_{\lbrack x_0-\mu,x_0+\mu \rbrack}(x), \hspace{0.3cm} x \in (0,\infty)
\end{eqnarray*}
where $s(\mu)$ is a scale to be determined later.
We denote the unscaled perturbation function by $\tilde{g}_\mu$; i.e.,
\begin{eqnarray*}
\tilde{g}_\mu(x) = (x_0+\mu -x)^{k+1}(x-x_0 + \mu)^{k+2}
1_{\lbrack x_0-\mu,x_0+\mu \rbrack}(x).
\end{eqnarray*}
If $\mu$ is chosen small enough so that the true density
$g_0$ is $k$-times differentiable on
$\lbrack x_0-\mu,x_0+\mu \rbrack $ and $g_0^{(k)}$ is
continuous on the latter interval, the perturbed function
$g_\mu$ is also $k$-times
differentiable on $\lbrack x_0-\mu, x_0+\mu \rbrack$
with a continuous $k$-th derivative.
Now, let $r$ be the function defined on $(0,\infty)$ by
\begin{eqnarray*}
r(x) =  (1-x)^{k+1}(1+x)^{k+2}1_{\lbrack -1,1 \rbrack}(x)
  =  (1-x^2)^{k+1}(1+x)1_{\lbrack -1,1 \rbrack}(x).
\end{eqnarray*}
Then, we can write $\tilde{g}_\mu$ as
\begin{eqnarray*}
\tilde{g}_\mu(x) = \mu^{2k+3}r\left( \frac{x-x_0}{\mu}\right).
\end{eqnarray*}
Then for $0 \leq j \leq k$
\begin{eqnarray*}
g^{(j)}_\mu(x_0) - g_0^{(j)}(x_0) = s(\mu) \mu^{2k+3-j}r^{(j)}(0).
\end{eqnarray*}
The scale $s(\mu)$ should  be chosen so that
for all $0 \leq j \leq k$
\begin{eqnarray*}
(-1)^j g^{(j)}_\mu(x) > 0, \hspace{0.3cm} \textrm{for}
\hspace{0.3cm} x \in \lbrack x_0-\mu, x_0+ \mu \rbrack.
\end{eqnarray*}
But for $\mu$ small enough, the sign of $(-1)^j g_\mu^{(j)}$
will be that of $(-1)^j g_0^{(j)}(x_0)$, and hence 
$g_{\mu}$ is $k-$monotone.  For $j=k$,
\begin{eqnarray*}
g^{(k)}_\mu(x_0)= g^{(k)}_0(x_0) + s(\mu) \mu^{k+3}r^{(k)}(0).
\end{eqnarray*}
Assume that $r^{(k)}(0) \ne 0$.
Set
\begin{eqnarray*}
s(\mu) = \frac{g_0^{(k)}(x_0)}{r^{(k)}(0)}\times \frac{1}{\mu^{k+3}} \, .
\end{eqnarray*}
Then for $ 0 \leq j  \leq k-1$
\begin{eqnarray*}
g_\mu^{(j)}(x_0)
 =  g_0^{(j)}(x_0)
        + \mu ^{k-j}\frac{g_0^{(k)}(x_0)r^{(j)}(0)}{r^{(k)}(0)} 
 =  g_0^{(j)}(x_0) + o(\mu),
\end{eqnarray*}
as $\mu \rightarrow 0$,
and so we can choose $\mu$ small enough so that
$(-1)^j g_\mu ^{(j)}(x_0) > 0 $. For $j=k$
\begin{eqnarray*}
(-1)^k g_\mu ^{(k)}(x_0) = 2 (-1)^k g_0^{(k)}(x_0) > 0.
\end{eqnarray*}
To show that $r^{(j)}(0) \ne 0$ for $0 \leq j \leq k$, we define
\begin{eqnarray*}
  x_{n,m} = \left((1-x^2)^{n} \right)^{(m)} \bigg \vert_{x=0}.
\end{eqnarray*}
Let $m \geq 2$ and $2n \geq m$. We have
\begin{eqnarray*}
\left((1-x^2)^{n} \right)^{(m)}
& = & \left(((1-x^2)^{n})'\right)^{(m-1)} \\
& = & \left(-2n x (1-x^2)^{n-1}  \right )^{(m-1)} \\
& = & -2n \left( x \left( (1-x^2)^{n-1}\right)^{(m-1)}
        + (m-1) \left((1-x^2)^{n-1} \right)^{(m-2)}  \right)
\end{eqnarray*}
where in the last equality, we used Leibniz's formula for
the  derivatives of a product; see e.g. \mycite{apost:57}, page 99.
Evaluating the last expression at $x=0$ yields
\begin{eqnarray*}
x_{n,m} = -2n (m-1) x_{n-1,m-2}.
\end{eqnarray*}
If $m$ is even, we obtain
\begin{eqnarray*}
x_{n,m}
& = & (-2)^{m/2} \prod_{i=0}^{m/2-1}(n-i) \times 
\prod_{i=0}^{m/2-1}(m-2i -1) \times
         x_{n-m/2,0} \\
& = & (-2)^{m/2} \prod_{i=0}^{m/2-1}(n-i) \times \prod_{i=0}^{m/2-1}(m-2i-1) \\
\end{eqnarray*}
since $x_{n-m/2,0}=1$.
Similarly, when $m$ is odd, we have
\begin{eqnarray*}
x_{n,m}
& = & (-2)^{(m-1)/2} \prod_{i=0}^{(m-1)/2-1}(n-i)
        \cdot \prod_{i=0}^{(m-1)/2-1}(m-2i-1)
\cdot x_{n-(m-1)/2,1} \\ & = & 0,
\end{eqnarray*}
since $x_{n-(m-1)/2,1}=0$.
Now, we have for $1 \leq j \leq k$
\begin{eqnarray*}
r^{(j)}(x)
& = & \left((1-x^2)^{k+1}(1+x)  \right)^{(j)} \\
& = & (x+1)\left( (1-x^2)^{k+1} \right)^{(j)}
         + j \left( (1-x^2)^{k+1}\right)^{(j-1)} \\
\end{eqnarray*}
and hence
\begin{eqnarray*}
r^{(j)}(0)
& = & \left((1-x^2)^{k+1} \right)^{(j)}_{x=0}
       + j \left((1-x^2)^{k+1}\right)^{(j-1)}_{x=0}.
\end{eqnarray*}
Therefore, when $j$ is even, the second term vanishes and
\begin{eqnarray*}
r^{(j)}(0) = (-2)^{j/2} \prod_{i=0}^{j/2-1}(k+1-i)
\times \prod_{i=0}^{j/2-1}(j-2i-1) \ne 0.
\end{eqnarray*}
When $j$ is odd, the first term vanishes and
\begin{eqnarray*}
r^{(j)}(0) & = & (-2)^{(j-1)/2}
     \prod_{i=0}^{(j-1)/2-1}(k+1-i) \times  j \times
       \prod_{i=0}^{(j-1)/2-1}(j-2i-2)\\
&= &  (-2)^{(j-1)/2} \prod_{i=0}^{(j-1)/2-1}(k+1-i)
       \times \prod_{i=0}^{(j-1)/2}(j-2i) \ne 0. \\
\end{eqnarray*}
We set
\begin{eqnarray*}
C_{k,j} = r^{(j)}(0), \hspace{0.3cm} \textrm{for} \hspace{0.3cm} 1 
\leq j \leq k \, .
\end{eqnarray*}
Then $C_{k,k}$ specializes to
\begin{eqnarray*}
C_{k,k} =\left \{
\begin{array}{ll}
(-2)^{k/2} \prod_{i=0}^{k/2-1}(k+1-i) \times
\prod_{i=0}^{k/2-1}(k-2i-1),
\hspace{1cm} \textrm{if $k$ is even}  \\
(-2)^{(k-1)/2}\prod_{i=0}^{(k-1)/2 -1}(k+1-i)
\times  \prod_{i=0}^{(k-1)/2}(k-2i),
\hspace{0.3cm} \textrm{if $k$ is odd}.
\end{array}
\right.
\end{eqnarray*}
The previous expressions can be given in a more compact form.
After some algebra, we find that
\begin{eqnarray}\label{FormulaCk}
C_{k,k} =\left \{
\begin{array}{ll}
  2 \times (-1)^{k/2} (k+1)(k-1)! {k \choose k/2-1},
\hspace{0.3cm}
\textrm{if $k$ is even}  \\
(-1)^{(k-1)/2} k! {k+1 \choose (k-1)/2},
\hspace{2.5cm} \textrm{if $k$ is odd}.
\end{array}
\right.
\end{eqnarray}
We have for $0 \leq j \leq k-1 $,
\begin{eqnarray*}
\vert T_j(g_\mu) - T_j(g_0) \vert
& = &  | g^{(j)}_{\mu}(x_0) - g_0^{(j)}(x_0) | \\
& = &  \left \vert \frac{C_{k,j}}{C_{k,k}}g_0^{(k)}(x_0) \right \vert  \mu^{k-j}
\equiv \lambda^{(j)}_{k,1}  \left \vert g_0^{(k)}(x_0) \right \vert  \mu^{k-j}
\end{eqnarray*}
where we defined
$\lambda^{(j)}_{k,1}
= \left \vert C_{k,j}/C_{k,k} \right \vert$ for $j \in \{ 0 , \ldots,  k-1 \}$.
Furthermore
\begin{eqnarray*}
\lefteqn{\int_0^\infty \frac{\left(g_\mu(x) - g_0(x)\right)^2}{g_0(x)}dx }\\
& = & \frac{\left(g_0^{(k)}(x_0)\right)^2}{\mu^{2(k+3)}
        (C_{k,k})^2} \int_{x_0-\mu}^{x_0 + \mu}
         \frac{(x_0 + \mu -x)^{2(k+1)}(x-x_0+\mu)^{2(k+2)}}{g_0(x)} dx \\
& = & \frac{\left(g_0^{(k)}(x_0)\right)^2}{\mu^{2(k+3)}(C_{k,k})^2}
        \int_{-\mu}^{\mu} \frac{(\mu^2 -y^2)^{2(k+1)}
         (y+\mu)^{2}}{g_0(x_0+y)} dy \\
& = & \frac{\left(g_0^{(k)}(x_0)\right)^2}{\mu^{2(k+3)}(C_{k,k})^2}
        \times \mu^{4(k+1)+3} \int_{-1}^{1}
         \frac{(1-z^2)^{2(k+1)}(z+1)^2}{g_0(x_0+\mu z)}dz \\
& = & \left(\frac{\left(g_0^{(k)}(x_0)\right)^2}{(C_{k,k})^2} \int_{-1}^{1}
        \frac{(1-z^2)^{2(k+1)}(z+1)^2}{g_0(x_0+\mu z)}dz \right) \mu^{2k+1} \\
& = & \left(\frac{\left(g_0^{(k)}(x_0)\right)^2}{g_0(x_0)} \frac{\int_{-1}^{1}
        (1-z^2)^{2(k+1)}(z+1)^2 dz}{(C_{k,k})^2} \right) \mu^{2k+1} + o(\mu^{2k+2})
\end{eqnarray*}
as $ \mu \searrow 0$.
This gives control of the Hellinger distance as well in view of
\mycite{jongb:00}, Lemma 2, page 282, or \mycite{jongb:95}, Corollary 3.2, pages
30 and 31.
We set
\begin{eqnarray*}
\lambda_{k,2} =  \frac{\int_{-1}^{1} (1-z^2)^{2(k+1)}(z+1)^2 dz}{(C_{k,k})^2}.
\end{eqnarray*}
The constants
$\lambda_{k,2} $ can be given more explicitly using the formula
\begin{eqnarray*}
I_{n,2p} = \int_0^1(1-x^2)^n x^{2p}dx
= 2^{2n+1}\frac{n!(n+1)!}{(2n+2)!}
\frac{{n+p \choose n+1}}{{2(n+p)+1 \choose 2(n+1)}},
\end{eqnarray*}
for any integers $n$ and $p$, using the convention
\begin{eqnarray*}
{n+p \choose n+1} = {2(n+p)+1 \choose 2(n+1)} = 1
\end{eqnarray*}
when $p=0$.
We have,
\begin{eqnarray*}
\int_{-1}^{1} (1-x^2)^{2(k+1)}(x+1)^2 dx
=  \int_{-1}^{1} (1-x^2)^{2(k+1)}x^2 dx + \int_{-1}^1(1-x^2)^{2(k+1)}dx,
\end{eqnarray*}
since
\begin{eqnarray*}
\int_{-1}^{1} (1-x^2)^{2(k+1)}x dx =0,
\end{eqnarray*}
and hence
\begin{eqnarray}\label{FormulaInt}
\lefteqn{\int_{-1}^{1} (1-x^2)^{2(k+1)}(x+1)^2 dx  =  2(I_{2(k+1),2}
         + I_{2(k+1),0})} \nonumber\\
& = & 2^{4k+6}
          \frac{(2(k+1))!(2k+3)!}{(4k+6)!}
          \frac{{2k+3 \choose 2k+3}}{{4k+7 \choose 4k+6}}
           + \frac{2^{4k+5}\left((2(k+1))!\right)^2}{(4k+5)!} \nonumber \\
& = & 2^{4k+5}\frac{\left((2(k+1))!\right)^2}{(4k+6)!} \left 
(\frac{2(2k+3)}{4k+7}
           + (4k+6) \right) \nonumber \\
& = & 2^{4k+5}\frac{\left((2(k+1))!\right)^2}{(4k+7)!} \left( (4k+6)
            + (4k+6)(4k+7)  \right) \nonumber \\
& = & 2^{4k+5}\frac{\left((2(k+1))!\right)^2}{(4k+7)!} (4k+6)(4k+8) 
\nonumber \\
& = & 2^{4(k+2)} (2k+3)(k+2)\frac{\left((2(k+1))!\right)^2}{(4k+7)!}.
\end{eqnarray}
Combining (\ref{FormulaCk})
and (\ref{FormulaInt}), we find that $\lambda_{k,2}$ is given by
\begin{eqnarray*}
\lambda_{k,2} =
2^{4(k+1)} \frac{(2k+3)(k+2)}{(k+1)^2}
\frac{\left((2(k+1))! \right)^2}{ (4k+7)! ((k-1)!)^2
\left( {k\choose k/2-1} \right )^2} ,
\hspace{0.3cm} \textrm{when $k$ is even,} \\
\end{eqnarray*}
and
\begin{eqnarray*}
\lambda_{k,2} = 2^{4(k+2)} (2k+3)(k+2)
\frac{\left((2(k+1))! \right)^2}{ (4k+7)!(k!)^2 
\left({k+1 \choose (k-1)/2} \right)^2},
\hspace{0.3cm} \textrm{when $k$ is odd.}
\end{eqnarray*}
\medskip
Now, by using the change of variable $\epsilon = \mu ^{2k+1} (b_{k} + o(1))$,
where
\begin{eqnarray*}
b_{k} = \lambda_{k,2} \frac{\bigg(g^{(k)}_0(x_{0})\bigg)^{2}}{g_0(x_{0})}
\end{eqnarray*}
so that
$\mu = \left(\epsilon/ b_k \right)^{1/(2k+1)}(1 + o(1))$,
then for $0 \leq j \leq k-1$, the modulus of continuity,
$m_j$, of the functional $T_j$ satisfies
\begin{eqnarray*}
m_{j}(\epsilon) \geq \lambda^{(j)}_{k,1}g^{(k)}_0(x_0)
\left(\frac{\epsilon}{b_k}\right)^{(k-j)/(2k+1)}(1 + o(1)).
\end{eqnarray*}
The result is that
\begin{eqnarray*}
m_{j}(\epsilon) \geq (r_{k,j} \epsilon )^{\frac{k-j}{2k+1}} (1+o(1)),
\end{eqnarray*}
where
\begin{eqnarray*}
r_{k,j} & = & \frac{\left(\lambda^{(j)}_{k,1}
g^{(k)}_0(x_0)\right)^{(2k+1)/(k-j)}}{b_k} \\
\end{eqnarray*}
and hence
\begin{eqnarray}
\lefteqn{\sup _{\tau > 0} \lim _{n \to \infty} \inf  n^\frac{k-j}{2k+1}
MMR_{1}(n,T_{j},\mathcal{D}_{k, n,\tau})} \nonumber  \\
& \geq & \frac{1}{4} \bigg(4 \frac{k-j}{2k+1}e^{-1} 
\bigg)^{\frac{k-j}{2k+1}}\left(r_{k,j}\right)^{\frac{k-j}{2k+1}},
\label{firstform}
\end{eqnarray}
which can be rewritten as
\begin{eqnarray*}
\lefteqn{\sup_{\tau > 0} \lim _{n \to \infty}
\inf n^{\frac{k-j}{2k+1}}MMR_{1}(n,T_{j},\mathcal{D}_{k, n,\tau})}\\
& \geq & \frac{1}{4}
       \bigg(4 \frac{k-j}{2k+1}e^{-1} \bigg)^{\frac{k-j}{2k+1}}
       \frac{\lambda^{(j)}_{k,1}}{\left(\lambda_{k,2}\right)^{\frac{k-j}{2k+1}}}
       \bigg \{ \left
\vert g^{(k)}_0(x_{0})\right 
\vert^{\frac{2j+1}{2k+1}}g_0(x_{0})^{\frac{k-j}{2k+1}} \bigg \}
\end{eqnarray*}
for $j=0, \cdots, k-1$.    \hfill $\blacksquare$

\bigskip

\section{Preliminary numerical results}
\label{NumIllustr}

 From the standard Exponential distribution $Exp(1)$ we simulated two samples of
respective sizes $n=100$ and $n=1000$. For any fixed
$k \geq 1$, the Exponential density is $k$-monotone.
Based on each sample, we computed the LSE and MLE for $k=3$ and $k=6$ in both the 
direct and inverse problems  using the iterative $(2k-1)$-th spline algorithm described 
in \mycite{balabwell:04b}. It should be noted that the true mixing distribution that
corresponds to a standard Exponential when viewed as a
$k$-monotone density is  $Gamma(k+1,1)$. Indeed,
\begin{eqnarray*}
\int_x^\infty \frac{1}{\Gamma(k)} (t-x)^{k-1} e^{-(t-x)} dt =1
\end{eqnarray*}
for all $x > 0$, and hence
\begin{eqnarray*}
\exp(-x)
& = & \int_x^\infty  \frac{(t-x)^{k-1}}{(k-1)!} e^{-t} dt
     =  \int_0^\infty  \frac{(t-x)^{k-1}_{+}}{(k-1)!} e^{-t} dt \nonumber \\
&= & \int_0^\infty  k \frac{(t-x)^{k-1}_{+}}{t^k} \frac{1}{k!} t^k e^{-t} dt
     =  \int_0^\infty  k \frac{(t-x)^{k-1}_{+}}{t^k} f_k(t)dt,
\end{eqnarray*}
where $f_k$ is the $Gamma(k+1,1)$ density. \\

\begin{figure}[H]
\hspace{0.0cm}
\epsfysize=5in 
\epsfxsize=5.5in
 {\includegraphics[bb=33 57 770 553,clip=false, width=12cm]{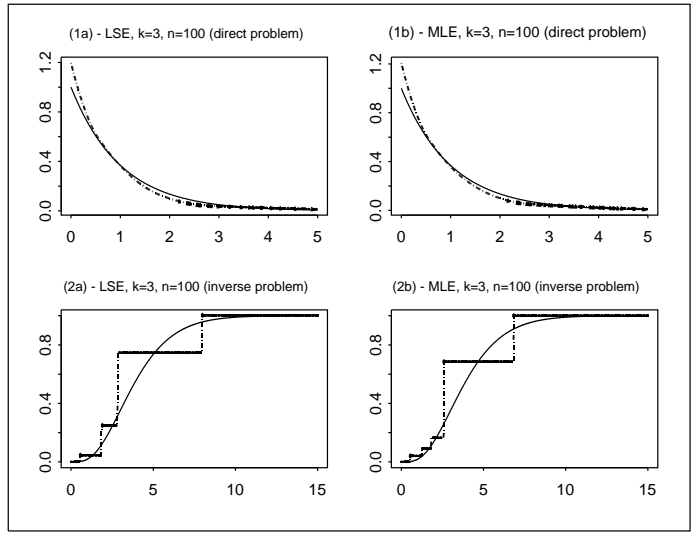}}
\caption{Illustration of $k$-montone estimation for $k=3$ via the ML and 
LS methods based on a  sample size $n=100$. Plots (1a) and (1b) show the 
LS and ML estimators (dashed lines)  of the exponential density (solid 
line). Plots (2a) and (2b) show the LS and ML estimators (dashed line) 
of $Gamma(4,1)$ (solid line), the true mixing distribution.}
\label{DirandInvK3n100}
\end{figure}

\begin{figure}[H]
\hspace{0.0cm}
\epsfysize=4in 
\epsfxsize=5.5in
{\includegraphics[bb=33 57 770 553,clip=true, width=12cm]{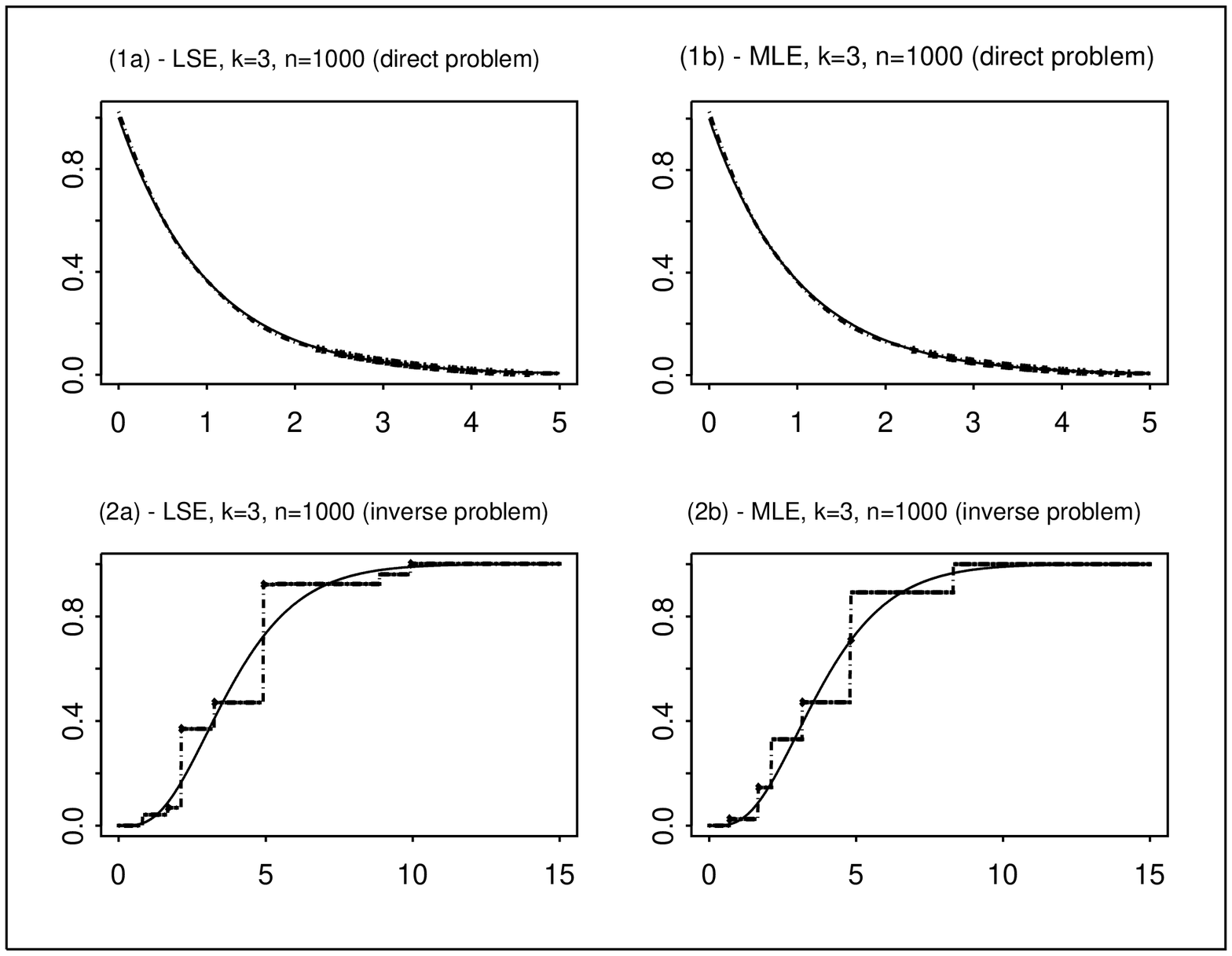}}
\caption{Illustration of $k$-montone estimation for $k=3$ via the ML and 
LS methods based on a sample size $n=1000$. Plots (1a) and (1b) show the 
LS and ML estimators (dashed lines)  of the exponential density (solid 
line). Plots (2a) and (2b) show the LS and ML estimators (dashed line) 
of $Gamma(4,1)$ (solid line), the true mixing distribution.}
\label{DirandInvK3n1000}
\end{figure}

\begin{figure}[H]
\hspace{0.0cm}
\epsfysize=4in 
\epsfxsize=5.5in
{\includegraphics[bb=33 57 770 553,clip=true, width=12cm]{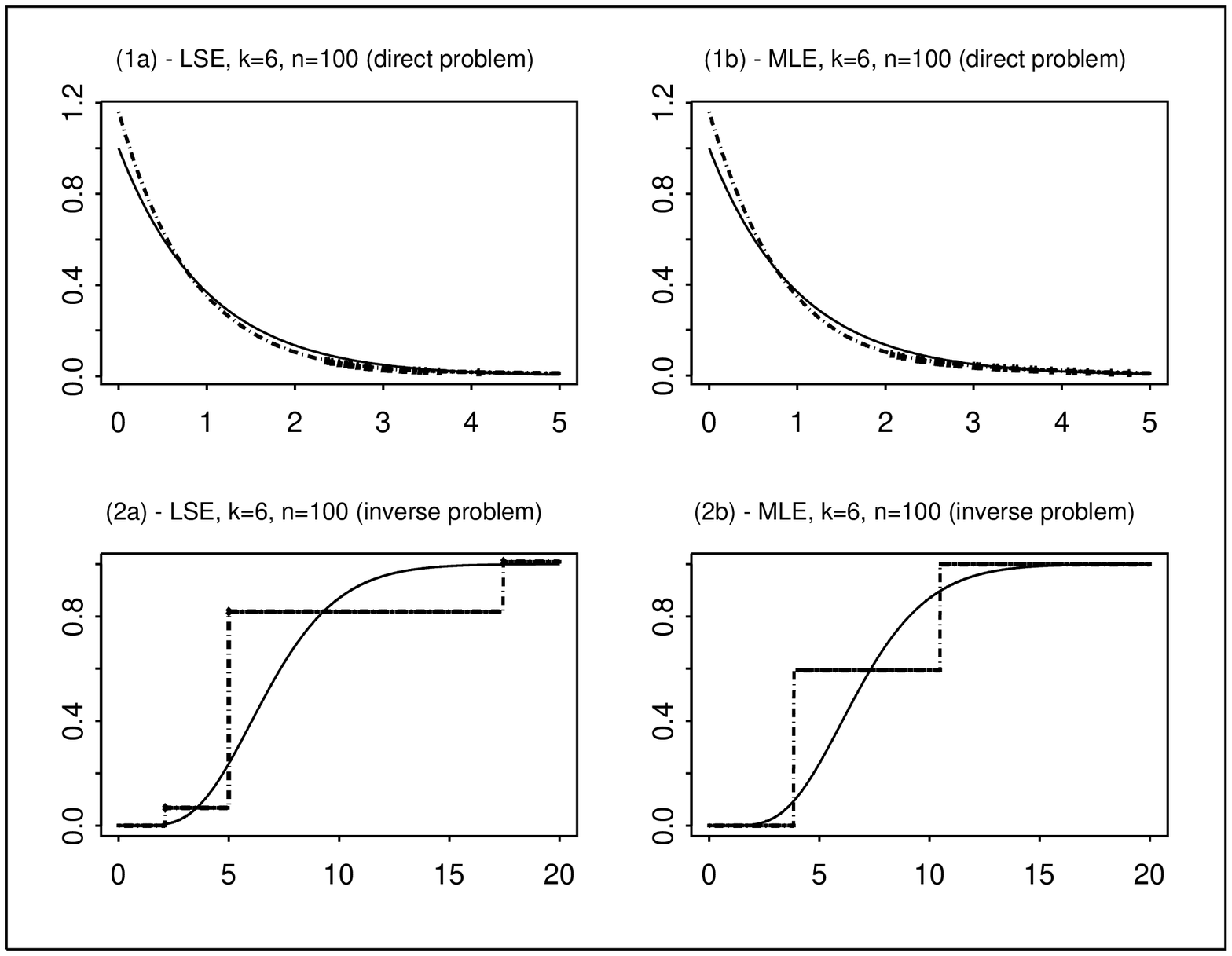}}
\caption{Illustration of $k$-montone estimation for $k=6$ via the ML and 
LS methods based on a sample size $n=100$. Plots (1a) and (1b) show the 
LS and ML estimators (dashed lines)  of the exponential density (solid 
line). Plots (2a) and (2b) show the LS and ML estimators (dashed line) 
of $Gamma(7,1)$ (solid line), the true mixing distribution.}
\label{DirandInvK6n100}
\end{figure}

\begin{figure}[H]
\hspace{0.0cm}
\epsfysize=4in 
\epsfxsize=5.5in
{\includegraphics[bb=33 57 770 553,clip=true, width=12cm]{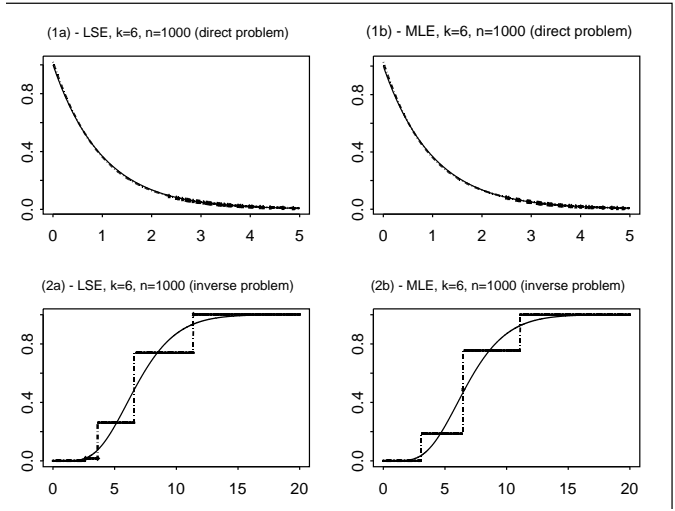}}
\caption{Illustration of $k$-montone estimation for $k=6$ via the ML and 
LS methods based on a sample size $n=1000$. Plots (1a) and (1b) show the 
LS and ML estimators (dashed lines)  of the exponential density (solid 
line). Plots (2a) and (2b) show the LS and ML estimators (dashed line) 
of $Gamma(7,1)$ (solid line), the true mixing distribution.}
\label{DirandInvK6n1000}
\end{figure}

For  $k=3$, the plots in Figures \ref{DirandInvK3n100} and 
\ref{DirandInvK3n1000} show the ML and LS estimators
of the Exponential density (direct problem) and the Gamma
distribution (inverse problem) based on $n=100$ and $1000$ respectively.
For $k=6$, similar plots were produced and are shown in
Figures \ref{DirandInvK6n100} and \ref{DirandInvK6n1000}.

\begin{table}[h]
\caption[Table of the LS estimates for $k=3, 6$
and $n=100, 1000$.]{Table of the obtained LS estimates
for $k=3, 6$ and $n=100, 1000$ and the corresponding
numbers of iterations $N_{it}$. A support point is denoted
by $\tilde{a}$ and its mass by $\tilde{w}$. }
\begin{center}
\begin{tabular}{ccc}
\hline
$k, n$ & $N_{it}$ & $(\tilde{a},\tilde{w})$   \\
\hline
$k=3, n=100$  & 13 &  $(0.569, 0.0459), (1.829, 0.168), (1.909,0.0347),$ \\
              &      &$(2.839,0.497), (7.939, 0.027), (7.989, 0.227)$ \\
$k=3, n=1000$ & 14 &   $(0.814,0.042), (1.674,0.027),(2.124, 0.300), (3.254,  0.100),  $ \\
              &  &  $(4.924,0.450), (5.334,0.001), (8.874, 0.037), (9.934, 0.039)$ \\
$k=6, n=100$ & 4 &  $(2.109, 0.067), (4.999, 0.750), (17.449, 0.190)$ \\
$k=6, n=1000$ & 6 &  $(2.625,0.017), (3.615,0.478), (6.575, 0.478), (11.375, 0.262)$  \\
\hline
\end{tabular}
\end{center}
\label{LSEresults}
\end{table}

\begin{table}[h]
\caption[Table of the ML estimates for
$k=3, 6$ and $n=100, 1000$.]{Table of the obtained
ML estimates for $k=3, 6$ and $n=100, 1000$. A support
point is denoted by $\hat{a}$ and its mass by $\hat{w}$. }
\begin{center}
\begin{tabular}{cc}
\hline
$k, n$ &  $(\hat{a},\hat{w})$   \\
\hline
$k=3, n=100$  &  $(0.549, 0.040), (1.259, 0.051), (1.819,0.072),$ \\
              &      $(2.579, 0.027), (2.589, 0.492), (6.839, 0.314)$ \\
$k=3, n=1000$  &   $(0.684,0.025), (1.664,0.120),(2.114, 0.184)$, \\
 &                    $(3.164,  0.141)$ \\
               &  $(4.794,0.236), (4.824,0.184), (8.304, 0.107)$ \\
$k=6, n=100$  &  $(3.839, 0.428), (3.849, 0.165), (10.479, 0.405)$ \\
$k=6, n=1000$ &  $(3.042, 0.186), (6.452, 0.300), (6.482, 0.267),$\\
              &     $ (11.072, 0.018), (11.102, 0.226)$  \\
\hline
\end{tabular}
\end{center}
\label{MLEresults}
\end{table}

\par\noindent

The figures illustrate consistency  in both the direct and
inverse problems, and it can be seen that convergence in the
direct problem is faster than it is in the inverse problem. 
This is already predicted by the corresponding theoretical 
rates of convergence, $n^{-k/(2k+1)}$ and $n^{-1/(2k+1)}$ respectively.

Note that the number of jump points of the estimators
of the mixing Gamma distribution, which are also the knots
of the estimators of the Exponential density, are fewer for
$k=6$ than for $k=3$: e.g. for $n=1000$, there are 8 jump
points for $k=3$ versus $4$ only when $k=6$ (for both estimators).
This was also observed in other simulations, and we
obtained even fewer points for larger values of $k$. This is
not surprising and is rather a consequence of the fact that
gap between the knots
(of  order $ n^{-1/(2k+1)} $)
is expected to get bigger with $k$. When $k$ increases, the
number of constraints on the estimated mixed
density grows, and hence it becomes harder
 to ``untangle'' the mixing distribution $F$ from the
very smooth Beta kernel. Finally, it should be mentioned
that although the MLE and LSE show very small visible
differences in the direct problem, it can be easily checked
by comparing the locations of jump points or the heights
of the jumps that these estimators are different
(compare Table \ref{LSEresults} and Table \ref{MLEresults}). \\

\medskip

\medskip

\section{Conclusion}
\label{conclusion}
In this first part, we have established existence of the MLE $\hat{g}_n$ and LSE $\tilde{g}_n$
of a $k$-monotone density $g_0$, and provided characterizations.
We have proved that both estimators are consistent in several senses  
as a first step toward understanding their asymptotic behavior. 
Consistency of higher derivatives of the estimators is usually not guaranteed in 
nonparametric density estimation problems, but here it is obtained 
``for free'' because of the particular shape constraints and 
smoothness of the density. 
In the sense of pointwise mean absolute error,  local asymptotic minimax 
lower bounds show that 
the rate of convergence of the $j$-th derivative of the MLE and LSE for 
$j=0, \cdots, k-1$ cannot be faster than $n^{-(k-j)/(2k+1)}$. 

Parts 3 and 4 are devoted to show that this rate, modulo a 
conjecture about boundedness of the error in a particular 
Hermite interpolation problem, is attained by the $j$-th derivative 
of the estimators, and that the joint asymptotic distribution of these 
derivatives involve a $(2k)$-convex stochastic process staying 
above (below) the $(k-1)$-fold integral of two-sided Brownian 
motion plus a deterministic drift if $k$ is even (odd). In the joint 
limiting distribution, the asymptotic variances are found to have 
the same dependence on $g_0(x_0)$ and $\vert g^{(k)}_0(x_0) \vert$ 
as the asymptotic constants obtained in the minimax lower bounds.
\bigskip

\par\noindent
{\bf Acknowledgements:}
We gratefully acknowledge helpful conversations with 
Carl de Boor, Nira Dyn, Tilmann Gneiting, and Piet Groeneboom.  
%

{}

\end{document}